\documentclass[final,12pt,3p]{elsarticle} 
\usepackage{import}
\usepackage{thesis}

\usepackage{graphicx}

\usepackage{amssymb}

\usepackage{lineno}

\journal{Engineering Structures}

\begin{document}

\pagenumbering{arabic}

\fancyfoot[C]{\thepage}

\begin{frontmatter}

\noindent \url{https://doi.org/10.1016/j.engstruct.2017.05.065}\\
This accepted manuscript is licensed under the CC-BY-NC-ND 4.0 license \faCreativeCommons\ \faCreativeCommonsBy\ \faCreativeCommonsNc\ \faCreativeCommonsNd\\

\noindent Fedorova, Maria, and M.V. Sivaselvan. 2017. “An Algorithm for Dynamic Vehicle-Track-Structure Interaction Analysis for High-Speed Trains.” \textit{Engineering Structures} 148: 857–77
\\

\title{An Algorithm for Dynamic Vehicle-Track-Structure Interaction Analysis for High-Speed Trains}

\author[label1]{Maria Fedorova}
\ead{mariafed@buffalo.edu}
\author[label1]{M.V. Sivaselvan\corref{cor1}}
\address[label1]{Department of Civil, Structural, and Environmental Engineering, University at Buffalo, Buffalo NY, USA}

\cortext[cor1]{Corresponding author}

\ead{mvs@buffalo.edu}



\begin{abstract}
The objective of the present work is to develop a robust, yet simple-to-implement algorithm for dynamic vehicle-track-structure-interaction (VTSI) analysis, applicable to trains passing over bridges. 
The algorithm can be readily implemented in existing bridge analysis software with minimal code modifications.
It is based on modeling the bridge and train separately, and coupling them together by means of kinematic constraints. The contact forces between the wheels and the track become Lagrange multipliers in this approach. A direct implementation of such an approach results in spurious oscillations in the contact forces. Two approaches are presented to mitigate these spurious oscillations - (a) a cubic B-spline interpolation of the kinematic constraints in time, and (b) an adaptation of an alternate time-integration scheme originally developed by Bathe. Solutions obtained using this algorithm are verified using a generic differential algebraic equation (DAE) solver. Due to high train speeds and possible track irregularities, wheels can momentarily lose contact with the track. This contact separation is formulated as a Linear Complementarity Problem (LCP). With this formulation, including contact separation in the analysis amounts to replacing a call to a linear equation solver by a call to an LCP solver, a modification of only two steps of the procedure. The focus of this paper is on the computational procedure of VTSI analysis. The main contribution of this paper is recognizing computational issues associated with time-varying kinematic constraints, clearly identifying their cause and developing remedies.

\end{abstract}


\end{frontmatter}


\section{Introduction} \label{sec:intro2D}

Vehicle-track-structure-interaction (VTSI) consists of the reciprocal influence of a bridge and a train on each other. As a train traverses a bridge, the deflections of the bridge as well as irregularities in the track act as support displacement input to the train at its wheels. The ensuing dynamics of the train in turn cause time-varying forces and vibration in the bridge. The purpose of VTSI analysis is to assure track safety and passenger comfort \cite{TM}. Passenger comfort is related to the acceleration experienced while the train is passing over the bridge~\cite{TM}. Track safety depends on the rate of loading of the track with time-varying forces from the train. For high-speed trains, this rate of loading may coincide with natural frequencies of the bridge resulting in resonance and amplification of bridge response. Moreover, geometric track irregularities may cause magnification of the contact forces and damage of the wheels and the track~\cite{neves2012direct}, or even train derailment under extreme conditions, such as earthquakes~\cite{Montenegro2016} and collisions~\cite{Xia2014}. Therefore, assessment of the vehicle running safety is another important reason for VTSI analysis.

Mathematically, VTSI can be represented as a system of two sets of equations of motion, for the train and bridge subsystems. Sun et al.~\cite{Sun} distinguish three types of algorithms for solving this system. The first group of algorithms aims to solve the system directly \cite{Xia2005,neves2012direct,Dimitrakopoulos2015,Salcher2015}. This approach is based on combining the two sets of equations into a single equation and solving the obtained equation. The second method requires condensation of vehicle degrees of freedom (DOF) into the bridge equation of motion and solving this updated equation. Based on such an approach, Yang et al.~\cite{Yang1999} proposed a vehicle-bridge interaction element. The third approach is to solve equations of motion of the vehicle and the bridge separately using iterative procedures \cite{Diana1989, Kwark2004, Zhang2013}. As Sun et al.~\cite{Sun} have pointed out, the first type of algorithms cannot be easily implemented into existing structural analysis software due to the fact that vehicle and bridge models are combined together. The second method does not allow to incorporate various train models into analysis and also requires specialized analysis software. The third approach is the most suitable in terms of incorporating a VTSI algorithm into existing software. However, as it was observed by Yang et al.~\cite{Yang1999}, the VTSI problem involves a large number of contact points between the wheels and track, hence convergence of iterative procedures may be low. Sivaselvan et al.~\cite{sivaselvan2014} proposed an algorithm that overcomes these drawbacks and can be integrated into existing software. The idea is to complement the system of equations of motion with a constraint equation and then solve the equations of motion of the vehicle and bridge separately using the constraint condition and avoiding an iterative procedure. In this case, a system of differential-algebraic equations (DAE) is obtained, that requires careful consideration of the numerical integration scheme. Moreover, adopting such a modular approach, contact separation between the wheels and the bridge can be easily modeled by formulating a linear complementary problem (LCP). The formulation of contact problems using complementarity methods has a long history~\cite{Ferris1997}, and, as a result, rolling wheel-rail contact can be modeled at various levels of detail ~\cite{Kwak1991,Xi2016_2D,Xi2016_3D,acary2008numerical}. Zhu et al.~\cite{Zhu2015} proposed a similar approach applying the mode superposition method to the bridge. 

In the present work, the VTSI algorithm proposed by Sivaselvan et al.~\cite{sivaselvan2014} is expanded upon. The main goal is to develop a highly modular algorithm that can be incorporated into existing software without interfering with the bridge model formulation. However, the cost of this modularity is the necessity to solve a DAE system. As opposed to the approach proposed by Zhu et al.~\cite{Zhu2015}, the current algorithm employs a general finite element model of the bridge directly, without utilizing the mode superposition method, which allows to broaden the usage of algorithm and apply it to different types of bridges. 

The organization of the paper is as follows. In Section~\ref{2D}, a system of governing equations of motion is derived. Various time-integration approaches are discussed in Sections~\ref{sec:Bauchau}-\ref{sec:Bathe}. Finally, in Section \ref{sec:separation}, contact separation between the wheel and the bridge is enabled.

\section{Governing equations} \label{2D}

\begin{figure}[h]
\centering
\includegraphics[scale=0.30]{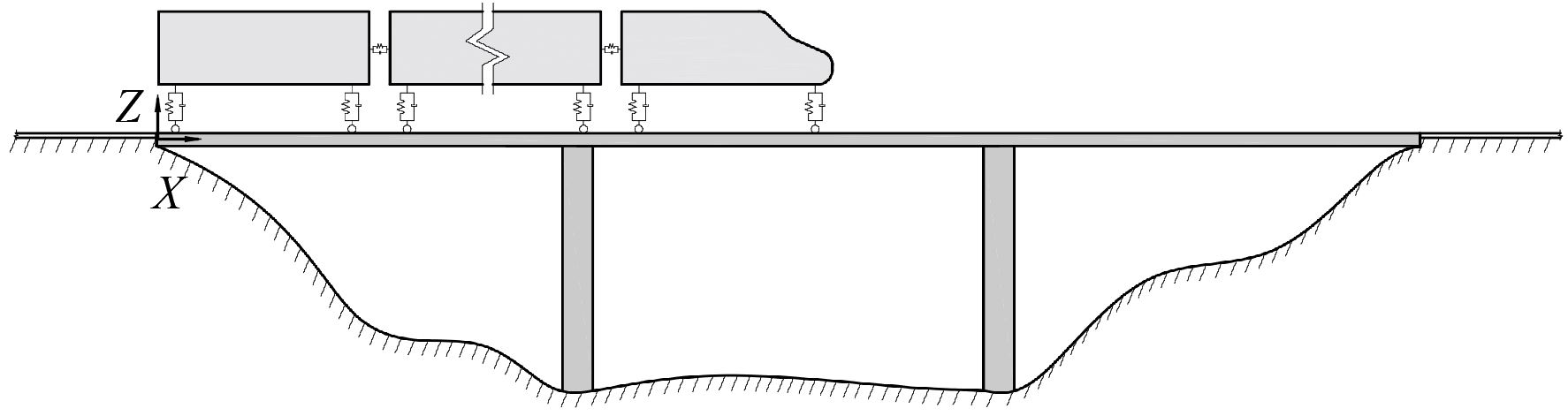}
\caption{Conceptual model of a train passing a bridge}
\label{fig:concept2D}
\end{figure}

A conceptual model of a train passing over a bridge is illustrated in Figure \ref{fig:concept2D}. The train is modeled as a sequence of cars. Each car is represented as a multibody system composed of rigid bodies, springs and dashpots. The bridge is modeled using standard structural or finite elements, such as beam and plate elements, box girders, cables, etc., that is any elements that can be utilized in commonly used bridge analysis software.

\begin{figure}[h]
\centering
\captionsetup{justification=centering}
\includegraphics[scale=0.3]{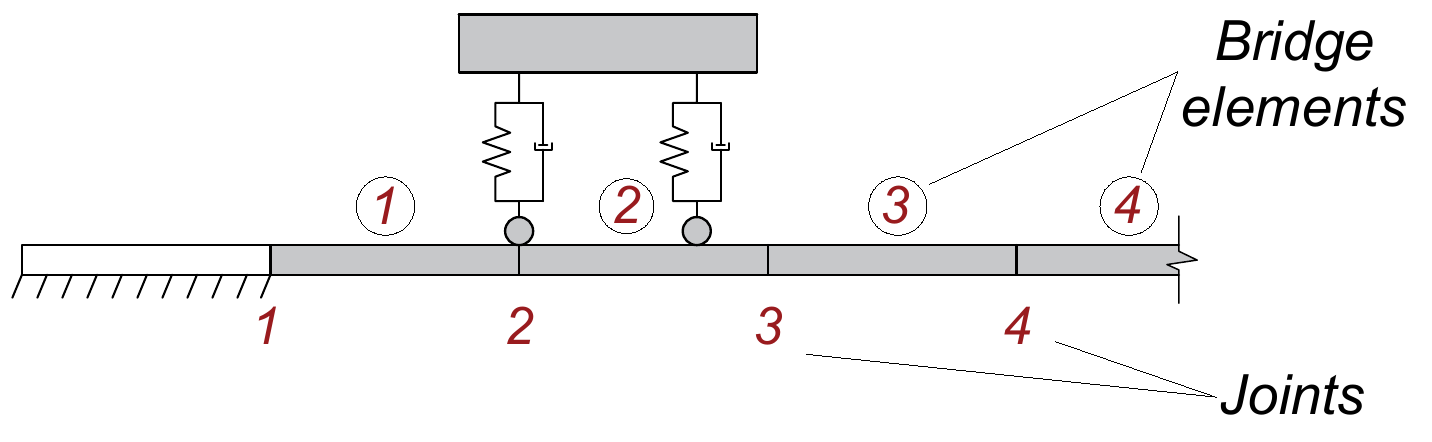}
\caption{One of the train wheels is on joint 2, hence it will be arbitrarily assigned to one of the two bridge elements connected to that joint, \circled{1} or \circled{2} }
\label{fig:joint_detail}
\end{figure}

The following assumptions are used to develop the algorithm for two-dimensional VTSI analysis \cite{sivaselvan2014}:

\begin{enumerate}
\item The train moves on a straight-line path along the bridge (curved paths are the topic of current research and will be addressed in a subsequent paper), so that the train dynamics is entirely planar, and the train loads on the bridge are in the global $Z$ direction.
\item The train and bridge displacements are small, so that linearized kinematics can be used for both. Material behavior is also linear.
\item The speed of the train is constant.
\item Through Section \ref{sec:Bathe}, the train wheels do not lose contact with the bridge (contact separation is considered in Section \ref{sec:separation}).
\item When a wheel is outside the span of the bridge, its displacement is zero.
\item At any time instant, each wheel is on only one element. If a wheel is on a joint, it is arbitrarily assigned to a bridge element connected to that joint (see Figure \ref{fig:joint_detail}).
\item Dead loads are applied to the bridge and the train using static analysis before the dynamic analysis is performed.
\end{enumerate}

\subsection{Train model} \label{sec:2D:train}

The train is modeled as a sequence of cars, each of which is a multibody system. The equation of motion of the train can then be written as \eqref{eq:train2D:short}. 
\begin{equation} \label{eq:train2D:short}
\Mt \utddotB + \Ct \utdotB + \Kt \utB + \left( \Lt \right)^T \lambdaB  = \Pt
\end{equation}
where $\Mt$, $\Ct$ and $\Kt$ are the mass, damping and stiffness matrices of the train model, $\Pt$ is the external load vector, such as a constant self-weight load on the train model, and $\utB$  is a vector of train displacements. Matrix $\Lt$ represents the influence of the reaction forces from the bridge on the train model. This matrix also plays a role in the constraint equation \eqref{eq:constraint2D}. Vector $\lambdaB$ is the vector of contact forces between the train wheels and the bridge (positive downward on the bridge and upward on the wheels). Here superscript ``t'' stands for train.

\subsubsection{Example}

\begin{figure}[h]
\centering
\captionsetup[subfigure]{justification=centering}
\begin{tabular}{l@{\hskip 1.0cm}r}
\subfloat[Parameters and degrees \newline of freedom; \newline ${\lct}$ is the total length of the car] 
{\label{fig:train2D:dof} \includegraphics[scale=0.32]{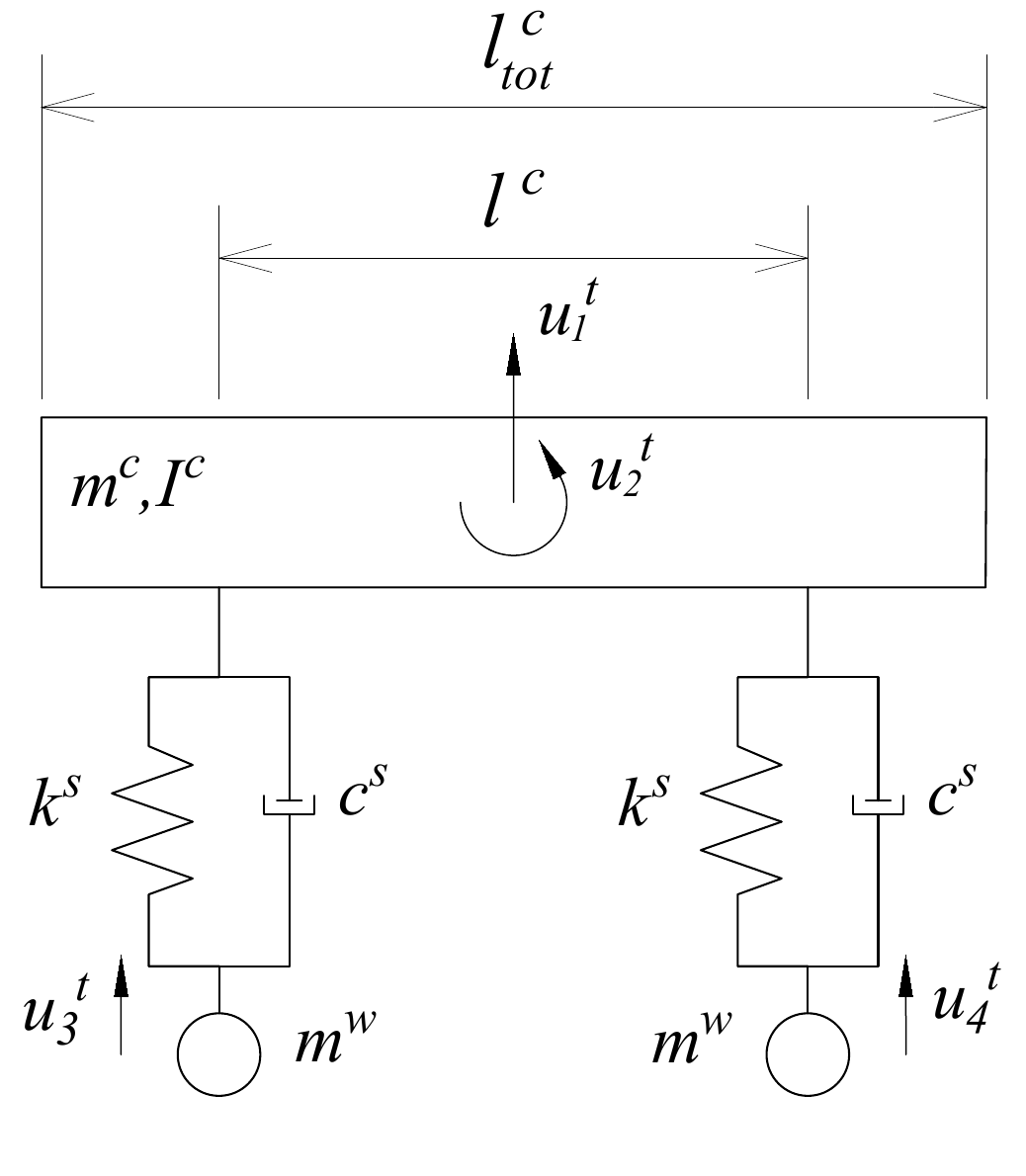}} &
\subfloat[Free body diagrams] 
{\label{fig:train2D:freebody} \includegraphics[scale=0.7]{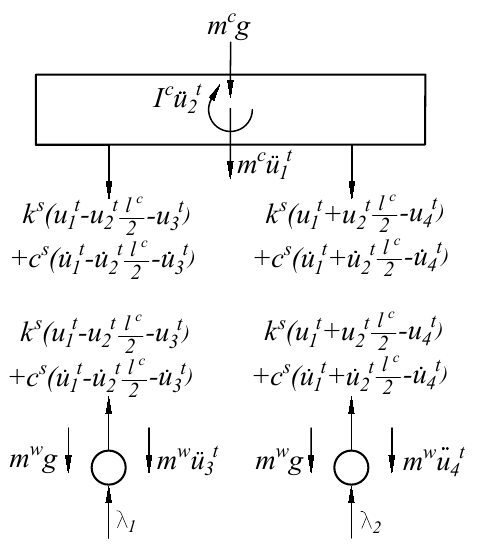}}
\end{tabular}
\caption{Simplified car model used for illustration}
\label{fig:train2D}
\end{figure}

To illustrate how equation \eqref{eq:train2D:short} comes about, a simple car model composed of a rigid bar connected to the wheels through dashpots and springs is considered. 
One such car is shown in Figure \ref{fig:train2D:dof}. The train is assumed to have \Ntwh\ wheels.

Using free body diagrams of the car and the wheels in Figure \ref{fig:train2D:freebody}, the equations of motion of the train model are derived as (compare with Equation \eqref{eq:train2D:short}):
\begin{equation} \label{eq:train2D}
\begin{gathered}                      
\left[ {\begin{array}{*{20}{c}}
{\mc}&{}&{}&{}\\
{}&{{\Ic}}&{}&{}\\
{}&{}&{\mw}&{}\\
{}&{}&{}&{{\mw}}
\end{array}} \right]\left( {\begin{array}{*{20}{c}}
{\utddot_1}\\
{\utddot_2}\\
{\utddot_3}\\
{\utddot_4}
\end{array}} \right) + \left[ {\begin{array}{*{20}{c}}
{2{\cs}}&0&{ - {\cs}}&{ - {\cs}}\\
0&{\frac{{{\cs}{\lc}^2}}{2}}&{\frac{{{\cs}{\lc}}}{2}}&{ - \frac{{{\cs}{\lc}}}{2}}\\
{ - {\cs}}&{\frac{{{\cs}{\lc}}}{2}}&{{\cs}}&0\\
{ - {\cs}}&{ - \frac{{{\cs}{\cs}}}{2}}&0&{{\cs}}
\end{array}} \right]\left( {\begin{array}{*{20}{c}}
{\utdot_1}\\
{\utdot_2}\\
{\utdot_3}\\
{\utdot_4}
\end{array}} \right)\\
 + \left[ {\begin{array}{*{20}{c}}
{2{\ks}}&0&{ - {\ks}}&{ - {\ks}}\\
0&{\frac{{{\ks}{\lc}^2}}{2}}&{\frac{{{\ks}{\lc}}}{2}}&{ - \frac{{{\ks}{\lc}}}{2}}\\
{ - {\ks}}&{\frac{{{\ks}{\lc}}}{2}}&{{\ks}}&0\\
{ - {\ks}}&{ - \frac{{{\ks}{\lc}}}{2}}&0&{{\ks}}
\end{array}} \right]\left( {\begin{array}{*{20}{c}}
{\ut_1}\\
{\ut_2}\\
{\ut_3}\\
{\ut_4}
\end{array}} \right) + \left[ {\begin{array}{*{20}{c}}
0&0\\
0&0\\
{ - 1}&0\\
0&{ - 1}
\end{array}} \right]\left( {\begin{array}{*{20}{c}}
{{\lambda _1}}\\
{{\lambda _2}}
\end{array}} \right) = \left( {\begin{array}{*{20}{c}}
{ - {\mc}g}\\
0\\
{ - {\mw}g}\\
{ - {\mw}g}
\end{array}} \right)\\
\end{gathered}
\end{equation}
where ${\mw}$ is the mass of each wheel; ${\mc}$ and ${\Ic}$ are the mass and moment of inertia of the carriage; $\ks$ and ${\cs}$ are the stiffness and damping of each suspension, ${\lc}$ is wheel-to-wheel distance; $\ut_i$ are displacements of the train model.

\subsection{Bridge model} \label{sec:2D:bridge}

The bridge is modeled employing any structural or finite elements available in commonly used bridge analysis software. The equation of motion of a bridge model can then be written in general form as
\begin{equation} \label{eq:bridge2D}
\Mb \ubddot + \Cb \ubdot + \Kb \ub = \Pb - {\left( {\Lb(t)} \right)^T}\lambdaB 
\end{equation}
where $\Mb$, $\Cb$ and $\Kb$ are the mass, damping and stiffness matrices of the bridge model; $\Pb$ is the external load vector, and $\utB$  is a vector of displacements of the bridge nodes. Matrix $\Lb$ is a time-dependent influence matrix derived based on the current locations of train wheels. Here superscript ``b'' stands for bridge.

\subsubsection{Example} \label{sec:2D:bridge:example}

\begin{figure}[h]
\centering
\includegraphics[scale=0.4]{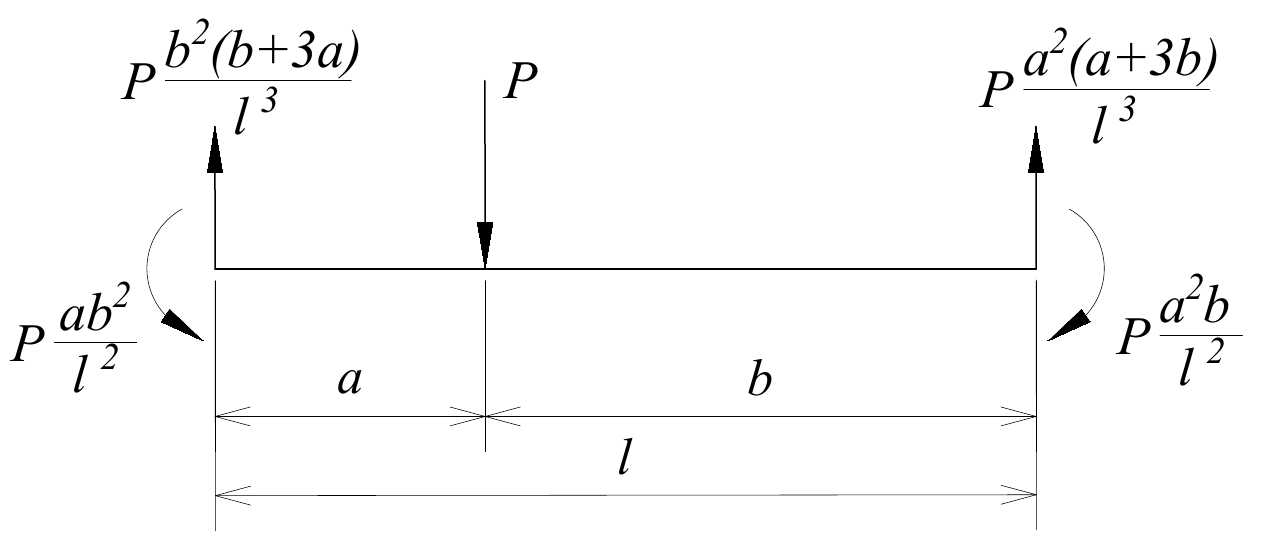}
\caption{Beam element with equivalent nodal forces}
\label{fig:bridge2D}
\end{figure}

An example of a bridge model assembled with prismatic beam elements is presented in this section for illustration. The bridge model is taken to have \Nbdof\ degrees of freedom. Figure \ref{fig:bridge2D} shows a typical beam element with equivalent nodal forces resulting from a point load P representing a train-wheel force. The equivalent nodal forces (fixed-end forces) for member $m$ can be written as
\begin{equation} 
{\bf F}_{\text{m}} = {\left( \Lbm \right)^T}P
\end{equation}
where
\begin{equation} \label{eq:Lb_mem_2D}
\Lbm = \left[ {\begin{array}{*{20}{c}}
{\frac{{{b^2}\left( {b + 3a} \right)}}{{{l^3}}}}&{\frac{{a{b^2}}}{{{l^2}}}}&{\frac{{{a^2}\left( {a + 3b} \right)}}{{{l^3}}}}&{ - \frac{{{a^2}b}}{{{l^2}}}}
\end{array}} \right]
\end{equation}
Vector \Lbm\ is the forces at the nodes due to a unit point load at the location ${(a,b)}$ in the element. Let ${\ubm}$ be the vector of member-end displacements. Then, by the principle of virtual work, ${\Lbm \ubm}$ is the displacement at the location ${(a,b)}$. 

A global influence matrix \Lb\ is then obtained by assembling element influence matrices \Lbm . For each wheel, the beam element on which the wheel is currently located is identified. The matrices \Lbm\ are calculated for the identified beam elements only, and placed in the global matrix \Lb\ in the row corresponding to the wheel and the columns corresponding to the degrees of freedom associated with the element. The matrix \Lb\ is therefore of size $\Ntwh \times \Nbdof$ and time-dependent. Then, by multiplying \Lb\ by the vector \ub\ of all the bridge node displacements, the bridge displacements at the current wheel location are obtained.

\subsubsection{Kinematic constraints and coupled equations of motion} \label{sec:2D:coupled}

The kinematic constraint can be written as
\begin{equation} \label{eq:constraint2D}
\Lt \utB + \Lb (t) \ub = -\rhoB \left( t \right)
\end{equation}
where ${\rhoB \left( t \right)}$ is the vector of elevation (vertical) track irregularities. Track irregularities amplify the dynamic impact of the wheel to the track. As a result, the wheel and the track may be damaged. For this reason, track irregularities should be incorporated into VTSI analysis.

\begin{figure}[h]
\centering
\begin{tabular}{l@{\hskip 0.5cm}r}
\subfloat[Contact~between~the~wheel~and~the~track] 
{\label{fig:constr_equal:real} {\makebox[.45\textwidth]{\includegraphics[scale=0.27]{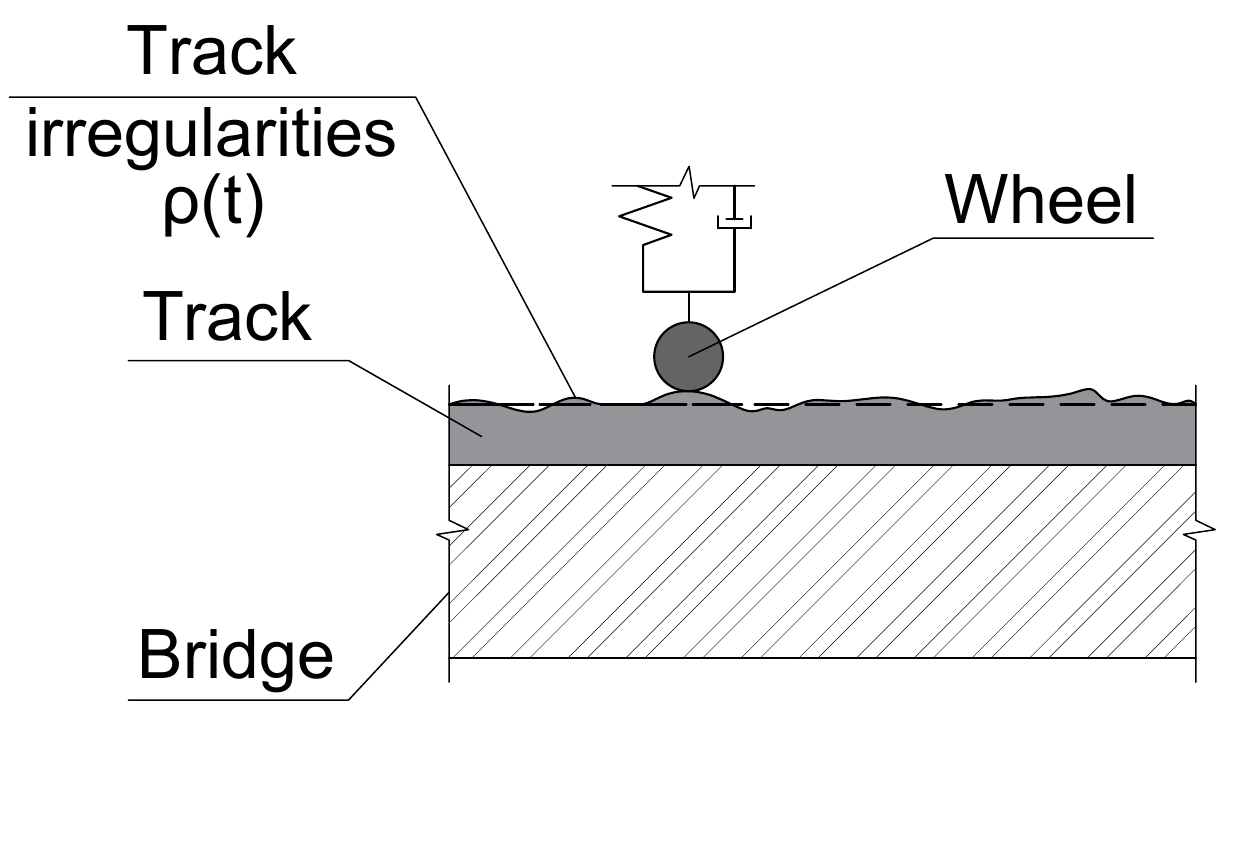}}}} &
\subfloat[Simplified modeling of the constraint] 
{\label{fig:constr_equal:model} \includegraphics[scale=0.3]{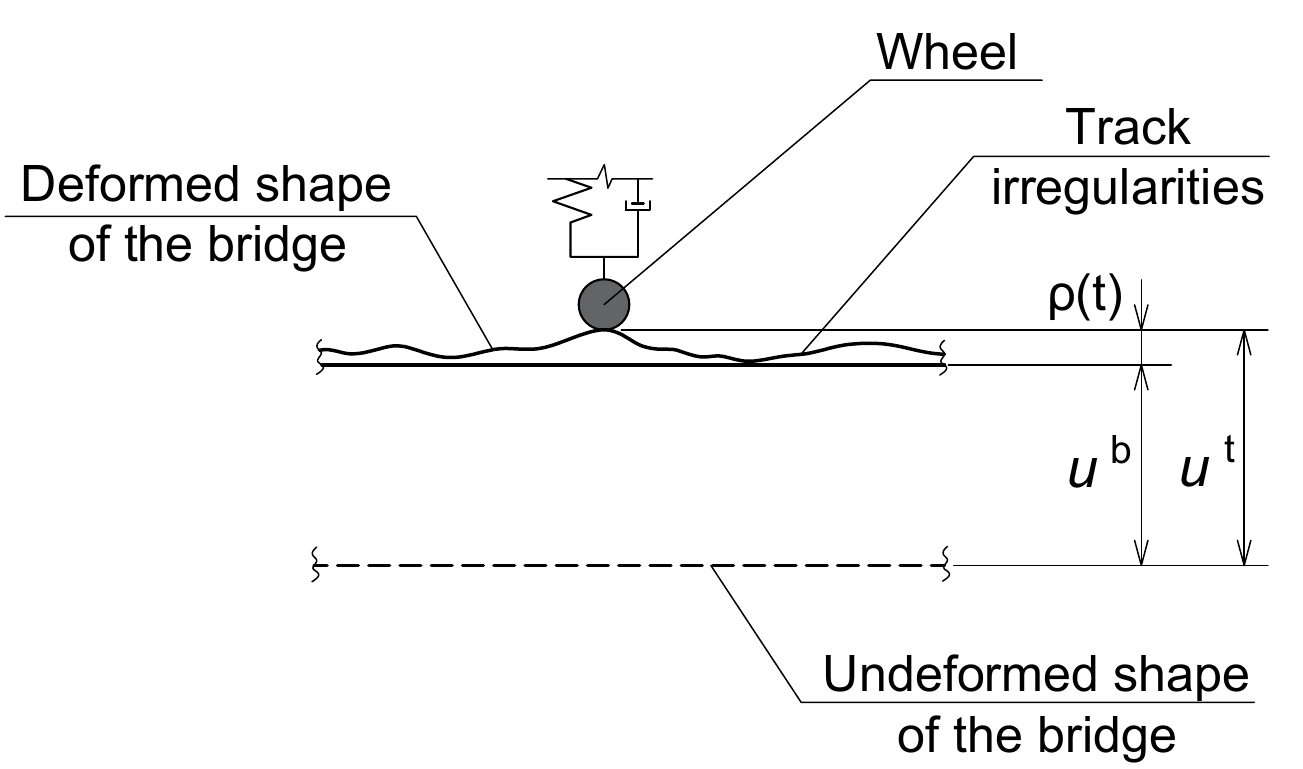}}
\end{tabular}
\caption{Modeling the kinematic constraint}
\label{fig:constr_equal}
\end{figure}

Figure \ref{fig:constr_equal:real} illustrates the contact between the wheel and the track. As can be seen, the wheel is constrained to the track. Since in this study the track is not modeled separately from the bridge, a simplified model of the constraint is adopted, as shown in Figure \ref{fig:constr_equal:model} (similarly to Yang et al. \cite{Yang1999})).

\FloatBarrier

The governing equations of motion of the train \eqref{eq:train2D:short} and the bridge \eqref{eq:bridge2D}, together with the constraint equation \eqref{eq:constraint2D}, represent the coupled equations of motion

\begin{subequations} \label{eq:coupled2D}
\begin{gather}  
\Mt \utddotB +\Ct \utdotB + \Kt \utB + {\left( \Lt \right)^T} \lambdaB  = \Pt \label{eq:coupled2Da}\\
\Mb \ubddot + \Cb \ubdot + \Kb \ub + {\left( \Lb(t) \right)}^T \lambdaB  = \Pb \label{eq:coupled2Db}\\
\Lt \utB  + \Lb(t) \ub = -\rhoB \left( t \right) \label{eq:coupled2Dc}
\end{gather} 
\end{subequations}
System \eqref{eq:coupled2D} is an index 3 system of differential-algebraic equations (DAE), which requires a careful consideration for time integration \cite{brenan1996numerical}.

\subsection{Remarks on coupling equations of motion directly without Lagrange multipliers} \label{sec:2D:remark}

Equations \eqref{eq:coupled2D} can be coupled directly, avoiding Lagrange multipliers. For this purpose, matrices in the train model should be partitioned as follows
\begin{equation}
    \utB = \begin{pmatrix} \uc \\ \uw \end{pmatrix}; \quad
    \Mt = \begin{bmatrix} \Mc & \\ & \Mw \end{bmatrix}; \quad
    \Ct = \begin{bmatrix} \Cc & \Ccw \\ \Cwc & \Cw \end{bmatrix}; \quad
    \Kt = \begin{bmatrix} \Kc & \Kcw \\ \Kwc & \Kw \end{bmatrix}; \quad
    \Pt = \begin{pmatrix} \Pc \\ \Pw \end{pmatrix}
\end{equation}

Equation \eqref{eq:coupled2D} can then be written as

\begin{subequations} \label{eq:CoupledEq:Direct}
    \begin{gather}
        \begin{alignedat}{7}
            &\Mc\dduc&&{}+{}\Cc\duc& &{}+{}\Ccw\duw&&{}+{}\Kc\uc& &{}+{}\Kcw\uw&&          &&{}= \Pc \\
            &\Mw\dduw&&{}+{}\Cwc\duc&&{}+{}\Cw\duw& &{}+{}\Kwc\uc&&{}+{}\Kw\uw& &{}-\lambdaB&&{}= \Pw
        \end{alignedat} \\
        \Mb\ubddot + \Cb\ubdot + \Kb\ub + \Lb(t)^\top\lambdaB = \Pb \\
        \uw = \Lb(t)\ub \label{eq:CoupledEq:Direct:c}
    \end{gather}
\end{subequations}

Differentiating equation \eqref{eq:CoupledEq:Direct:c} twice, we get
\begin{equation}
    \begin{aligned}
        \duw  &= \dLb\ub + \Lb\ubdot \\
        \dduw &= \ddLb\ub + 2\dLb\ubdot + \Lb\ubddot
    \end{aligned}
\end{equation}

Eliminating \lambdaB\ and \uw\ in equations \eqref{eq:CoupledEq:Direct}, we get

\begin{subequations} \label{eq:CoupledEq:Direct_2}
    \begin{align}
        \text{Train car} &:\ \Mc\dduc + \Cc\duc + \Kc\uc = \Pc - \underbrace{(\Ccw\dLb+\Kcw\Lb)\ub}_{\text{\small{bridge displacement}}} - \underbrace{\Ccw\Lb\ubdot}_{\text{\small{bridge velocity}}} \label{eq:CoupledEq:Direct_2:a}\\
        \text{Bridge}         &:\ \begin{gathered}
                                    \del{\Mb+(\Lb)^\top\Mw\Lb}\ubddot + \del{\Cb+(\Lb)^\top\Cw\Lb + 2(\Lb)^\top\Mw\dLb}\ubdot \\
                                    + \del{\Kb + (\Lb)^\top\Kw\Lb + (\Lb)^\top\Cw\dLb + (\Lb)^\top\Mw\ddLb}\ub \\
                                    = \Pb + (\Lb)^\top\underbrace{(\Pw - \Cwc\duc - \Kwc\uc)}_{\text{\small{train force}}}
                                 \end{gathered}  \label{eq:CoupledEq:Direct_2:b}
    \end{align}
\end{subequations}

As can be seen from equation \eqref{eq:CoupledEq:Direct_2}, the bridge response influences the the train response, and at the same time, the train force influences the bridge response. Coupling the equations this way leads to a reduction of the size of the system and an elimination of the necessity to solve DAE. However, the bridge mass, damping and stiffness matrices are modified (compare equations \eqref{eq:coupled2Db} and \eqref{eq:CoupledEq:Direct_2:b}) and must be updated at each time step. Furthermore,  due to the time-dependence of the constraint between train and bridge displacements, these modified matrices are not symmetric. The second derivative \ddLb\ is also commonly not continuous when a wheel crosses a node from one element to another, because curvature is not necessarily continuous across beam elements. Moreover, as was described in Section \ref{sec:intro2D}, an algorithm based on such an approach cannot be easily implemented into existing structural analysis software due to the fact that vehicle and bridge models are combined together and modularity is lost. Therefore, in subsequent sections, the bridge and train are coupled together using Lagrange multipliers, per equations \eqref{eq:coupled2D}.
\section{Discretization using Bauchau scheme} \label{sec:Bauchau}

Equations (\ref{eq:coupled2D}) are discretized in time using a self-stabilized algorithm proposed by Bauchau \cite{bauchau2003} for constrained multibody systems. This scheme preserves the energy of the system, and therefore is unconditionally stable. Bauchau pointed out that enforcing position and velocity constraints simultaneously may over-constrain a system. Instead, he suggested imposing the governing equations of motion (\ref{eq:coupled2Da}) and (\ref{eq:coupled2Db}) at the half of time step ${n+1/2}$, and the kinematic constraint (\ref{eq:coupled2Dc}) at the end of time step ${n+1}$. According to this scheme, the velocity-displacement relationships are discretized as follow 
\begin{equation}
\begin{gathered}
\utB_{n + 1} = \utB_n + \Delta t\left( \frac{{\utdotB_{n + 1} + \utdotB_n}}{2} \right) \\
\ub_{n + 1} = \ub_n + \Delta t\left( \frac{{\ubdot_{n + 1} + \ubdot_n}}{2} \right)
\end{gathered}
\end{equation}
This discretization reduces to Newmark's method with parameters 0.25, 0.5 in the absence of constraints. It also should be noted that accelerations \utddotB\ and \ubddot\ are obtained not directly from the discretization scheme, but rather as a byproduct of computations at the half of time step ${n+1/2}$: $\utddotB_{n+1/2}$ and $\ubddot_{n+1/2}$ (see Procedure \ref{procedure_2D}).

The equations \eqref{eq:coupled2D} are then discretized as
\begin{equation} \label{discr2D}
\begin{gathered}
\Mt \left( \frac{\utdotB_{n + 1} - \utdotB_n}{\Delta t} \right) + \Ct \left( \frac{\utdotB_{n + 1} + \utdotB_n}{2} \right) + \Kt \left( \frac{\utB_{n + 1} + \utB_n}{2} \right) + {\Lt}^T \lambdaB_{n + 1/2} = \frac{\Pt_{n + 1} + \Pt_n}{2}\\
\Mb \left( \frac{\ubdot_{n + 1} - \ubdot_n}{\Delta t} \right) + \Cb \left( \frac{\ubdot_{n + 1} + \ubdot_n}{2} \right) + \Kb \left( \frac{\ub_{n + 1} + \ub_n}{2} \right) + {\Lb (t_{n+1/2})}^T \lambdaB_{n + 1/2} = \frac{\Pb_{n + 1} + \Pb_n}{2}\\
\Lt \utB_{n + 1} + \Lb \left( t_{n + 1} \right) \ub_{n + 1} = -\rhoB \left( t_{n + 1} \right)
\end{gathered}
\end{equation}

Velocity at the mid-point ${n+1/2}$ can be defined as
\begin{equation}
\begin{gathered}
\vtbar = \frac{\utdotB_{n + 1} + \utdotB_n}{2}\\
\vbbar = \frac{\ubdot_{n + 1} + \ubdot_n}{2}
\end{gathered}
\end{equation}

Discretized equations \eqref{discr2D} can then be rearranged into the following linear system of equations
\begin{equation} \label{discr2D_linear}
\begin{gathered}
\Mtbar \vtbar + \frac{\Delta t}{2} {\Lt}^T \lambdaB_{n + 1/2} = \at \\
\Mbbar \vbbar + \frac{\Delta t}{2} {\Lb (t_{n+1/2})}^T \lambdaB_{n + 1/2} = \ab \\
\frac{\Delta t}{2} \Lt \vtbar + \frac{\Delta t}{2} \Lb (t_{n+1}) \vbbar = {\bf b}
\end{gathered}
\end{equation}
where

\begin{subequations} \label{rearranged}
\begin{gather}  
\Mtbar = \Mt + \frac{\Delta t}{2} \Ct + \frac{{\Delta t}^2}{4}\Kt \label{rearranged_a} \\
\Mbbar = \Mb + \frac{\Delta t}{2} \Cb + \frac{{\Delta t}^2}{4}\Kb \label{rearranged_b} \\
\at = \Mt \utdotB_n + \frac{\Delta t}{2} \left( \frac{\Pt_{n + 1} + \Pt_n}{2} - \Kt \utB_n \right) \label{rearranged_c} \\
\ab = \Mb \ubdot_n + \frac{\Delta t}{2} \left( \frac{\Pb_{n + 1} + \Pb_n}{2} - \Kb \ub_n \right) \label{rearranged_d} \\
{\bf b} =  - \frac{1}{2} \left( \Lt \utB_n + \Lb (t_{n+1}) \ub_n \right) - \frac{1}{2} \rhoB \left( t_{n + 1} \right) \label{rearranged_e}
\end{gather} 
\end{subequations}

The system (\ref{discr2D_linear}) is solved for ${\left( \vtbar, \vbbar, \lambdaB_{n + 1/2} \right)}$ at each time step. The train model is computed separately and incorporated into a structural analysis code by means of calling a function once.
A Schur complement-type approach is used resulting in a coupling equation of size equal to the number of wheels (step~17 of Procedure~1). The vector of wheels positions $\xw(t)$ along the bridge (i.e. along the ${x}$ axis) computed based on the train speed $v^t$ and initial positions of the wheels $\xw(0)$ :
\begin{equation} \label{wheelpos}
\xw(t) = \xw(0) + v^t t
\end{equation}

The bridge matrices \Mb, \Cb, \Kb\ and \Mbbar\ are not time-dependent and therefore computed only once. The procedure of the analysis is outlined below.

\floatname{algorithm}{Procedure} 
\begin{algorithm}[H]
\caption{\text Linear dynamic analysis of bridge-train interaction using Bauchau scheme} \label{procedure_2D}
\begin{algorithmic}[1]
\State Given: \Pt, \Pb, \rhoB(t), ${\Delta t}$ 
\State Assemble \Mt, \Ct, \Kt and \Mtbar
\State Assemble \Mb, \Cb, \Kb and \Mbbar
\State Compute $\utB_0$ by linear static analysis. Set $\utdotB_0 = 0$ \Comment{Train initial conditions}
\State Compute $\ub_0$ by linear static analysis. Set $\ubdot_0 = 0$ \Comment{Bridge initial conditions}
\For{$n\leftarrow 0$, NUMINC}
\State Compute $\xw(t_{n+1/2})$ by equation (\ref{wheelpos})
\State Obtain $\Lb(t_{n+1/2})$
\State Compute $\xw(t_{n+1})$ by equation (\ref{wheelpos})
\State Obtain $\Lb(t_{n+1})$
\State Compute \at\ by equation \eqref{rearranged_c}
\State Compute \ab\ by equation \eqref{rearranged_d}
\State Solve $\Mtbar \left[ \vttilde\ \At \right] = \left[ \at\ {\Lt}^T \right]$
\State Solve $\Mbbar \left[ \vbtilde\ \Ab \right] = \left[ \ab\ {\Lb (t_{n+1/2})}^T \right]$
\State Set ${\bf A} = \dfrac{\dt^2}{2}\left( \Lt \At + \Lb \left( t_{n + 1} \right) \Ab \right)$
\State Set ${\boldsymbol {\bar \rho}}  = \Lt \left( \utB_n + \dt \vttilde \right) + \Lb \left( t_{n + 1} \right) \left( \ub_n + \dt \vbtilde \right) + \rhoB \left( t_{n + 1/2} \right)$
\State Solve ${\bf A} \lambdaB_{n + 1/2} = {\boldsymbol {\bar \rho}}$
\State Set $\vtbar = \vttilde - \frac{\dt}{2} \At \lambdaB_{n + 1/2}$; $\utdotB_{n + 1}= 2\vtbar - \utdotB_n$; $\utB_{n + 1} = \utB_n + \dt \vtbar$;  
\Statex \hspace*{1.2cm} $\utddotB_{n + 1/2} = \del{ \utdotB_{n+1} - \utdotB_n} / \dt$
\State Set $\vbbar = \vbtilde - \frac{\dt}{2} \Ab \lambdaB_{n + 1/2}$; $\ubdot_{n + 1}= 2\vbbar - \ubdot_n$; $\ub_{n + 1} = \ub_n + \dt \vbbar$; 
\Statex \hspace*{1.2cm} $\ubddot_{n + 1/2} = \del{ \ubdot_{n+1} - \ubdot_n} / \dt$
\EndFor
\end{algorithmic}
\end{algorithm}

\FloatBarrier
\subsection{Numerical example} \label{sec:Bauchau:ex}
In order to test this algorithm and determine possible shortcomings, an example of dynamic interaction between a two-span continuous bridge and one car is presented in this section. 

The material and geometric properties of the bridge are adopted from \cite{yang2004vehicle} and as follows: length of one span 30m, Young's modulus ${E=29}$GPa, cross-section second moment of inertia ${I_y=8.65 {\text m}^4}$, mass per unit length ${36000{\text{kg/m}}}$. Track irregularities are not considered, therefore ${\rhoB(t)=0}$. A Rayleigh damping matrix \Cb\ is constructed based on the damping ratio $\xi = 5\%$ in the first and second modes of the bridge. First two natural frequencies of the bridge model are equal to $7.2$Hz and $10.44$Hz for fixed bridge (Cases 1 and 2), and $4.61$Hz and $7.2$Hz for simply supported bridge (Cases 3 and 4).

The parameters of the train model are adopted from \cite{spiryagin2014design}. The train consists of one car with mass ${\mc=60000}$kg, wheel-to-wheel length ${\lc=15}$m, $\lct = 20$m, moment of inertia ${\Ic=1125{\text t}\cdot {\text m}^2}$. Each suspension has the stiffness ${\ks=5000{\text{kN/m}}}$. Damping of the system is assumed to be equal to $5\%$ of the critical damping. For each suspension damping \cs\ is equal to ${27{\text{kN}} \cdot {\text{ s/m}}}$. Natural frequency of the car is equal to $2.06$Hz. According to the Technical Memorandum 2.10.10 \cite{TM}, the maximum speed for which dynamic VTSI analysis should consider is 1.2 times the line design speed or $250$mph, whichever is less. Therefore, we assume that the speed of the train \vt\ is constant and equal to ${110{\text{m/s}}}$ ($246$mph). Influence of train speed is considered in Section \ref{sec:resonance}.

\begin{figure}[H]
\centering
\begin{tabular}{c}
\subfloat[Cases 1-4] 
{\label{fig:cases:1-4}\includegraphics[scale=0.58]{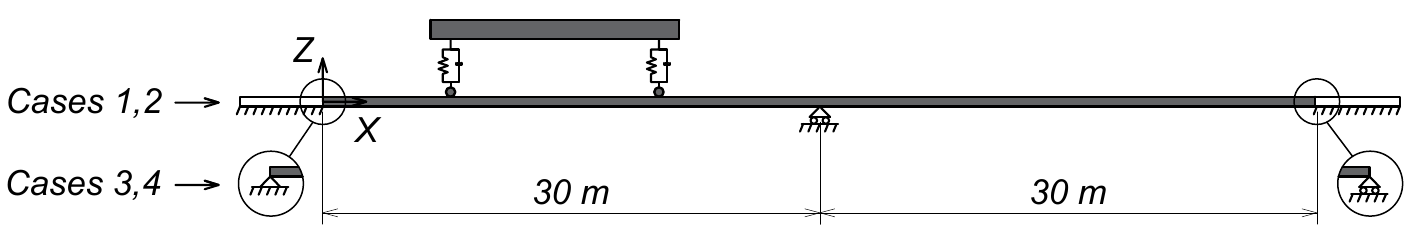}}\\
\subfloat[Case 5] 
{\hspace{0.8cm}{\label{fig:cases:5}\includegraphics[scale=0.50]{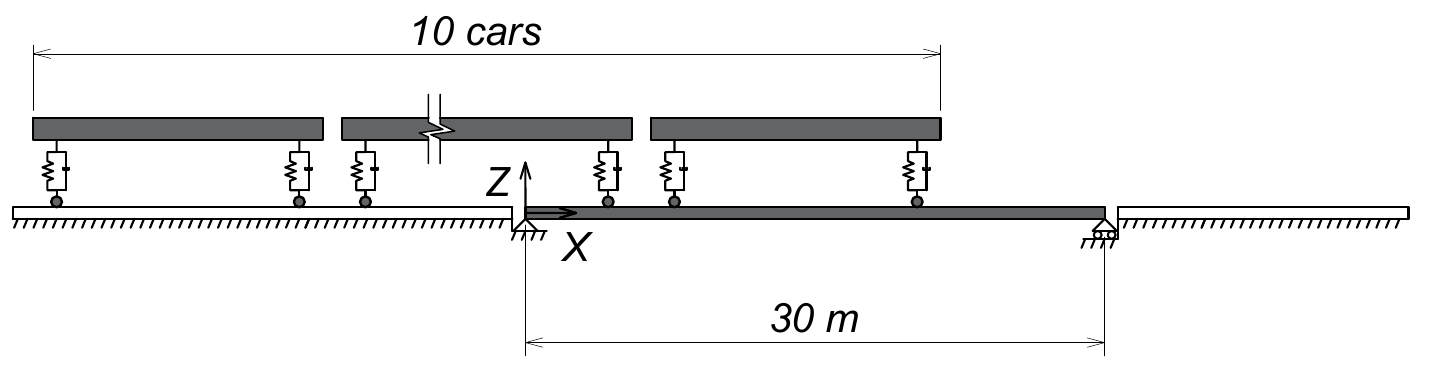}}}\\
\subfloat[Case 6] 
{\hspace{0.8cm}{\label{fig:cases:6}\includegraphics[scale=0.42]{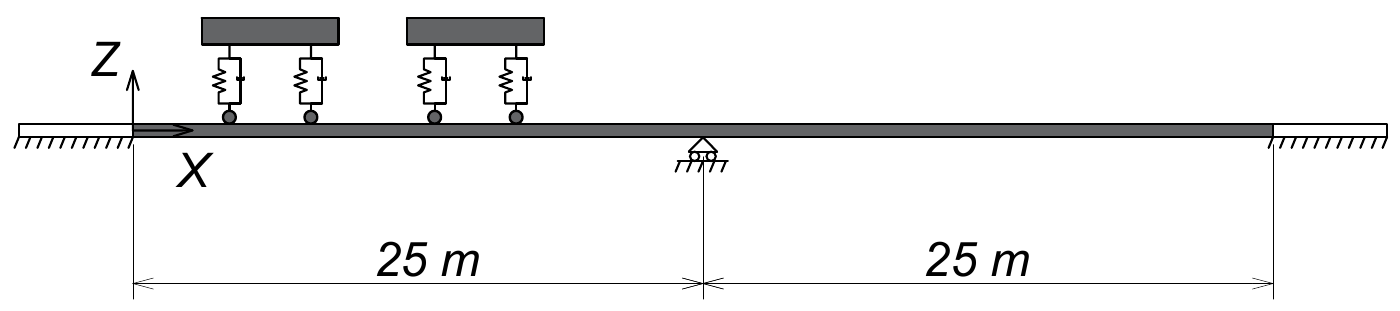}}} 
\end{tabular}
\caption{Numerical example: A bridge and a passing train}
\label{fig:cases}
\end{figure}

Four cases with different boundary conditions and different wheel masses \mw\ are examined (Cases 1 to 4 in the Table \ref{table:cases}) in this section. For each type of support condition, fixed and simply supported, two separate sub-cases are analyzed: mass of the wheel \mw\ is equal to $0$kg and ${\mw=1000}$kg. While the second value is realistic, the first value is introduced to explore a possibility of oscillations in the results. The numerical example was preliminarily evaluated with various masses of the wheel (e.g. 0, 100, 1000kg). Obtained results showed that spurious oscillations occur only for the mass of the wheel greater than 0kg, no matter how small the mass is. Therefore, the results for two representative masses (0 and 1000kg) are presented to estimate numerical issues. Cases 1-4 are illustrated in Figure \ref{fig:cases:1-4}. For cases 1 and 2, additional fixed elements are modeled at each end of the bridge in order to enable B-spline interpolation in Section \ref{sec:BSpline}. Cases 5 and 6 are presented in Figures \ref{fig:cases:5} and \ref{fig:cases:6} and described in detail in Sections \ref{sec:resonance} and \ref{sec:LCP:ex} correspondingly. 

\begin{table}[h!]
\centering
\begin{tabular}{ | m{1.0cm} | m{2.0cm} | m{0.9cm} | m{1.8cm} | m{1.8cm} | m{1.8cm} | m{1.8cm}| m{1.8cm} | } 
\hline
Cases & Support conditions & $\mw$, kg 
& Bauchau scheme \footnotesize{(Sec. \ref{sec:Bauchau:ex})} 
& B-spline \footnotesize{(Sec. \ref{sec:BSpline:ex})} 
& Bathe method \footnotesize{(Sec. \ref{sec:Bathe:ex})}
& LCP \hspace{0.5cm} \footnotesize{(Sec. \ref{sec:LCP:ex})}
& Sundials \hspace{0.5cm} \footnotesize{(Sec. \ref{sec:Sundials})}\\ 
\hline \hline
1 & Fixed & 0 & \checkmark & \checkmark & \checkmark & &   \\ 
\hline
2 & Fixed & 1000 & \checkmark & \checkmark & \checkmark & & \checkmark \\ 
\hline 
3 & Simply supported & 0 & \checkmark &  & \checkmark & &  \\ 
\hline 
4 & Simply supported & 1000 & \checkmark &  & \checkmark & &  \\ 
\hline
5 & Simply supported & 1000 &  &  & \checkmark & &  \\ 
\hline 
6 & Fixed & 1000 & & & & \checkmark &  \\ 
\hline
\end{tabular}
\caption{Analysis cases}
\label{table:cases}
\end{table}

\FloatBarrier

\begin{figure}[H]
\centering
\setcounter{subfigure}{0}
\captionsetup[subfigure]{justification=centering}
\begin{tabular}{l@{\hskip 1.5cm}r}
\subfloat[Maximum contact force] 
{\label{fig:converg:step:tr} \includegraphics{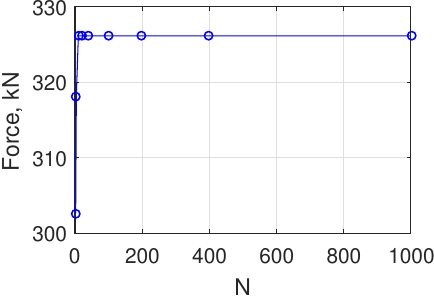}} &
\subfloat[Maximum contact force] {\label{fig:converg:time:tr} \includegraphics{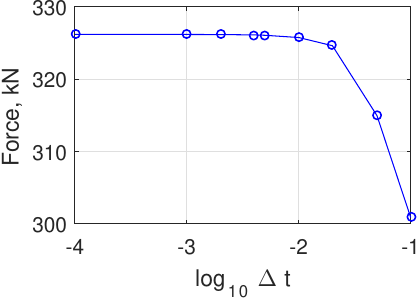}}\\[0pt]
\subfloat[Maximum displ. of the midpoint of the first span] 
{\label{fig:converg:step:br} \includegraphics{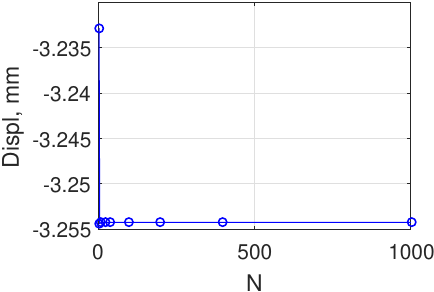}} & 
\subfloat[Maximum displ. of the midpoint of the first span] {\label{fig:converg:time:br} \includegraphics{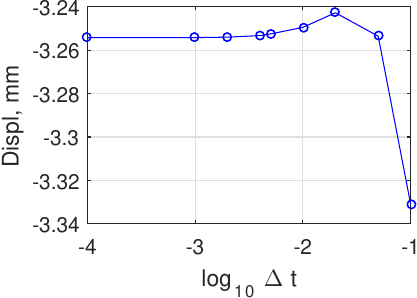}}\\[0pt]
\subfloat[Maximum displ. of the first wheel] 
{\label{fig:converg:step:force} \includegraphics{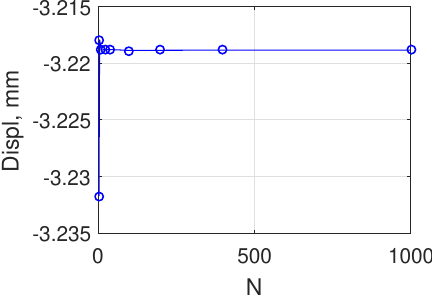}}&
\subfloat[Maximum displ. of the first wheel] {\label{fig:converg:time:force} \includegraphics{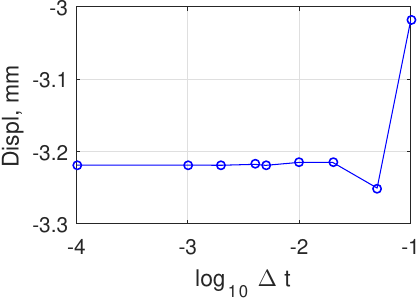}}\\
\end{tabular}
\caption{Convergence of the solution depending on number of elements $N$ and time step $\Delta t$}
\label{fig:converg}
\end{figure}

In order to explore the convergence of the solution, Case 1 is analyzed using different numbers of elements $N =$ 2, 4, 10, 20, 40, 100, 200, 400, 1000 to model each span of the bridge. Based on the obtained results (Figures \ref{fig:converg:step:tr}, \ref{fig:converg:step:br} and \ref{fig:converg:step:force}), each span is discretized with $N$ = 100 elements. Case 1 is then analyzed again considering different sizes of time step $\Delta{t}$ = 0.1, 0.05, 0.02, 0.01, 0.005, 0.004, 0.002, 0.001, 0.0001 s. As can be seen from Figures \ref{fig:converg:time:tr}, \ref{fig:converg:time:br} and \ref{fig:converg:time:force}, all three estimated parameters converge for time step $\Delta{t} =$ 0.001s. Therefore, $\Delta{t}$ is set to be equal to 0.001s for detailed analysis. 

\FloatBarrier
\textit{Cases 1 and 2}

\begin{figure}[H]
\centering
\setcounter{subfigure}{0}
\captionsetup{justification=centering}
\captionsetup[subfigure]{justification=centering}
\begin{tabular}{l@{\hskip 0.8cm}r}
\subfloat[Case 1, $\mw = 0$] 
{\includegraphics[scale =0.96]
{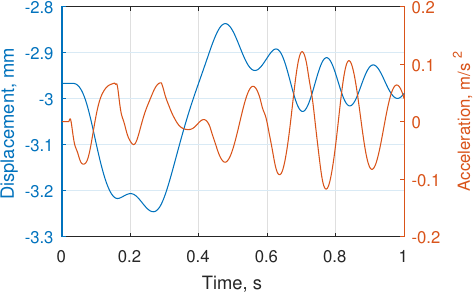}} &
\subfloat[Case 2, $\mw = 1000\ kg$] 
{\label{fig:Bauchau:1:br:2}\includegraphics[scale =0.96]
{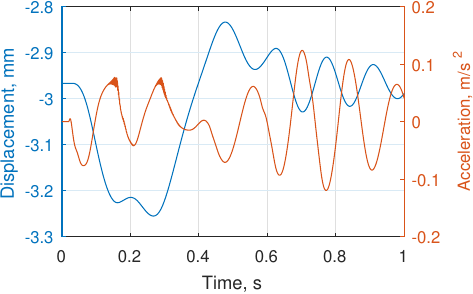}}
\end{tabular}
\caption{Bauchau scheme. Cases 1 (a) and 2 (b). Vertical displacements and accelerations of the middle point of the first span}
\label{fig:Bauchau:1:br}
\end{figure}

\begin{figure}[H]
\centering
\captionsetup{justification=centering}
\captionsetup[subfigure]{justification=centering}
\begin{tabular}{l@{\hskip 0.8cm}r}
\subfloat[Case 1, $\mw = 0$] 
{\includegraphics{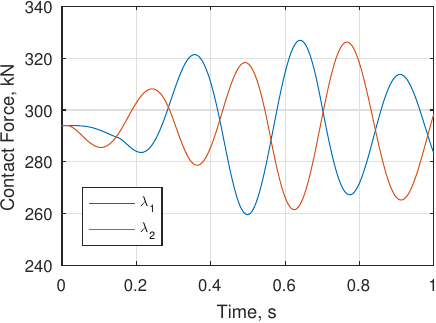}} &
\subfloat[Case 2, $\mw = 1000\ kg$] 
{\label{fig:Bauchau:1:forces:2}\includegraphics{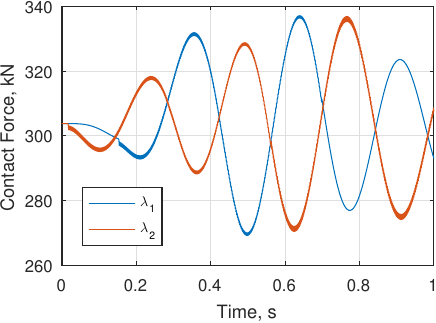}}
\end{tabular}
\caption{Bauchau scheme. Cases 1 (a) and 2 (b). Contact forces}
\label{fig:Bauchau:1:forces}
\end{figure}

\begin{figure}[H]
\centering
\captionsetup{justification=centering}
\captionsetup[subfigure]{justification=centering}
\begin{tabular}{l@{\hskip 0.8cm}r}
\subfloat[Case 2, $\mw = 1000\ kg$. \hspace{40mm} Displacements of the wheels] 
{\includegraphics{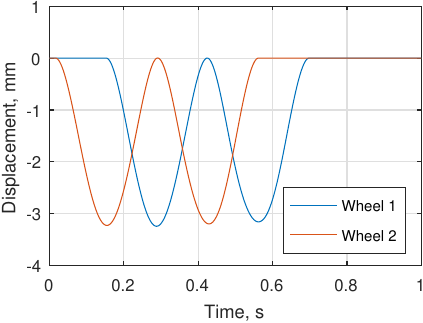}} &
\subfloat[Case 2, $\mw = 1000\ kg$. \hspace{40mm} Displacements and rotations of the car] 
{\includegraphics{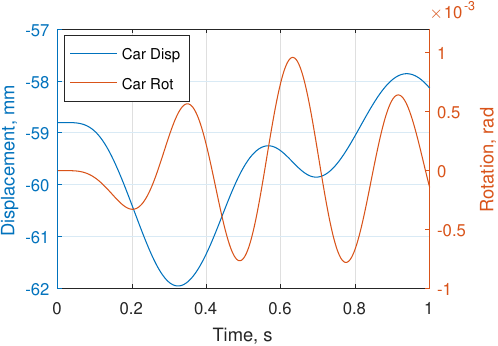}}
\end{tabular}
\caption{Bauchau scheme. Case 2. Vertical displacements of the car and the wheels, rotations of the car}
\label{fig:Bauchau:1:tr:displ}
\end{figure}

\begin{figure}[H]
\centering
\captionsetup{justification=centering}
\captionsetup[subfigure]{justification=centering}
\begin{tabular}{l@{\hskip 0.8cm}r}
\subfloat[Case 2, $\mw = 1000\ kg$. \hspace{40mm} Accelerations of the wheels] 
{\label{fig:Bauchau:1:tr:acc:wh}\includegraphics{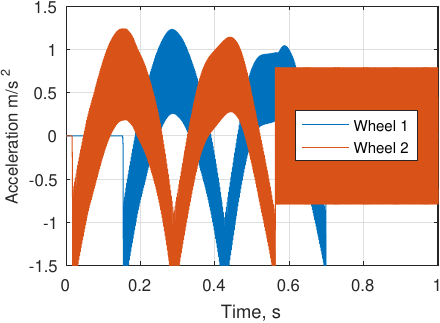}} &
\subfloat[Case 2, $\mw = 1000\ kg$. \hspace{40mm} Accelerations of the car] 
{\label{fig:Bauchau:1:tr:acc:2}\includegraphics{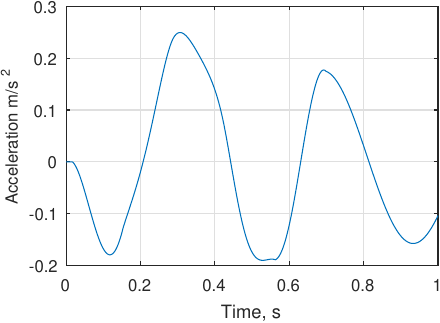}}
\end{tabular}
\caption{Bauchau scheme. Case 2. Accelerations of the wheels and car}
\label{fig:Bauchau:1:tr:acc}
\end{figure}

To illustrate vertical displacements and accelerations of the bridge, responses of the middle point of the first span are shown in Figure \ref{fig:Bauchau:1:br}. As can be seen, some spurious oscillations occur in bridge accelerations in Case~2. Contact forces between the wheels and the bridge are presented in Figure \ref{fig:Bauchau:1:forces}. For the case $\mw=0$ (Case 1), when the wheels play a role only in constraining the displacements, the implemented scheme performs well. However, the increase of \mw\ causes spurious oscillations in contact forces. Figures \ref{fig:Bauchau:1:tr:displ} and \ref{fig:Bauchau:1:tr:acc} show vertical displacements and accelerations of the train model. Both figures represent the analysis of the Case 2 ($\mw = 1000$kg). Plots of the displacements for Case 1 are not presented here, since the differences in the results are negligible because of the high stiffness of the bridge. Figure \ref{fig:Bauchau:1:tr:acc:wh} shows that spurious oscillations in the contact forces are significant for Case~2. For Case~1, no oscillations in accelerations were observed.

\FloatBarrier
\textit{Cases 3 and 4}

\begin{figure}[H]
\centering
\captionsetup{justification=centering}
\captionsetup[subfigure]{justification=centering}
\begin{tabular}{l@{\hskip 0.8cm}r}
\subfloat[Case 3, $\mw = 0$] 
{\includegraphics[scale =0.96]
{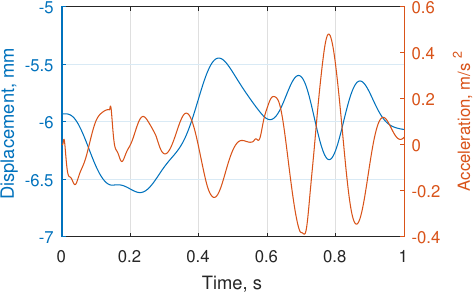}} &
\subfloat[Case 4, $\mw = 1000\ kg$] 
{\label{fig:Bauchau:0:br:4}\includegraphics[scale =0.96]
{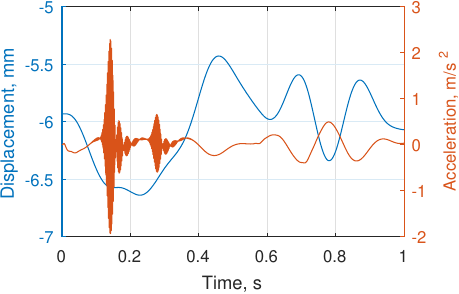}}
\end{tabular}
\caption{Bauchau scheme. Cases 3 (a) and 4 (b). Vertical displacements and accelerations of the middle point of the first span}
\label{fig:Bauchau:0:br}
\end{figure}

\begin{figure}[H]
\centering
\captionsetup{justification=centering}
\captionsetup[subfigure]{justification=centering}
\begin{tabular}{l@{\hskip 0.8cm}r}
\subfloat[Case 3, $\mw = 0$] 
{\includegraphics{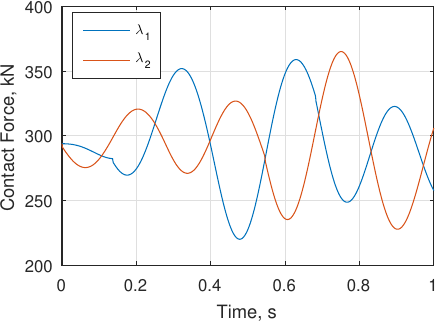}} &
\subfloat[Case 4, $\mw = 1000\ kg$] 
{\label{fig:Bauchau:0:forces:4}\includegraphics{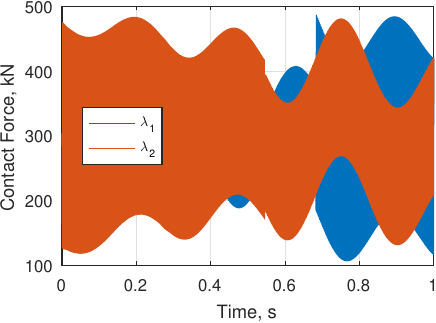}}
\end{tabular}
\caption{Bauchau scheme. Cases 3 (a) and 4 (b). Contact forces}
\label{fig:Bauchau:0:forces}
\end{figure}

\begin{figure}[H]
\centering
\captionsetup{justification=centering}
\captionsetup[subfigure]{justification=centering}
\begin{tabular}{l@{\hskip 0.8cm}r}
\subfloat[Case 4, $\mw = 1000\ kg$. \hspace{40mm} Displacements of the wheels] 
{\includegraphics{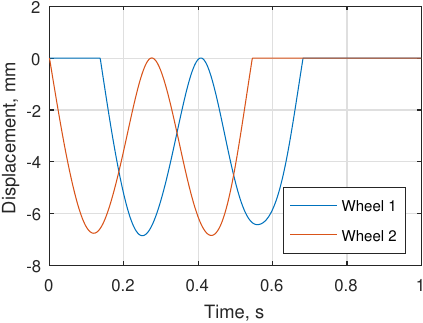}} &
\subfloat[Case 4, $\mw = 1000\ kg$. \hspace{40mm} Displacements and rotations of the car] 
{\includegraphics{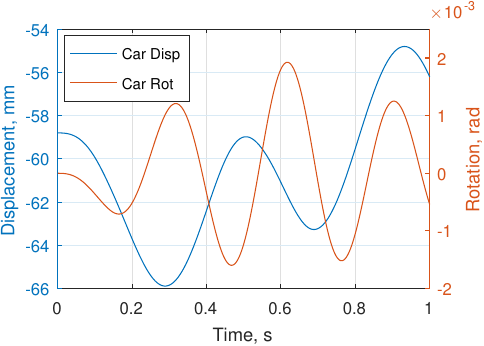}}
\end{tabular}
\caption{Bauchau scheme. Case 4. Vertical displacements of the car and the wheels, rotations of the car}
\label{fig:Bauchau:0:tr:displ}
\end{figure}

\begin{figure}[H]
\centering
\captionsetup{justification=centering}
\captionsetup[subfigure]{justification=centering}
\begin{tabular}{l@{\hskip 0.8cm}r}
\subfloat[Case 4, $\mw = 1000\ kg$. \hspace{40mm} Accelerations of the wheels] 
{\includegraphics{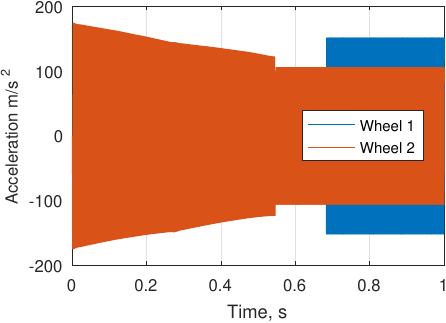}} &
\subfloat[Case 4, $\mw = 1000\ kg$. \hspace{40mm} Accelerations of the car] 
{\label{fig:Bauchau:0:tr:acc:4}\includegraphics{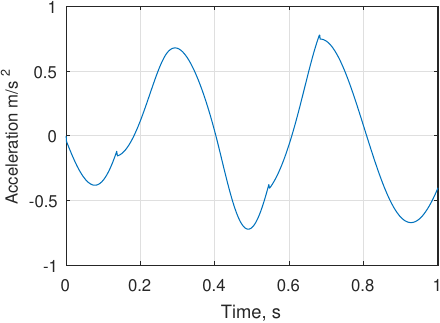}}
\end{tabular}
\caption{Bauchau scheme. Case 4. Acceleration of the wheels and car}
\label{fig:Bauchau:0:tr:acc}
\end{figure}

The results of analyses for Cases 3 and 4 are presented in Figures \ref{fig:Bauchau:0:br} - \ref{fig:Bauchau:0:tr:acc}. It can be seen that spurious oscillations in accelerations and contact forces are more significant than in the analogous Cases 1 and 2 with fixed support conditions (compare Figures \ref{fig:Bauchau:1:br:2} and \ref{fig:Bauchau:0:br:4}, \ref{fig:Bauchau:1:forces:2} and \ref{fig:Bauchau:0:forces:4}, \ref{fig:Bauchau:1:tr:acc:2} and \ref{fig:Bauchau:0:tr:acc:4}). These oscillations are numerical artifacts and have to be eliminated from the results of VTSI analysis in order to properly interpret the results and estimate wheel/track safety (based on the contact forces values) and passenger comfort (based on the train accelerations), and to validate the obtained results by experimental data.

For the current application, the second derivative of \Lb\ should be continuous, since the second derivative of Equation \eqref{eq:coupled2Dc} relates accelerations of the train and the bridge. However, the influence matrix \Lb\ is constructed using Hermite polynomials which provide only $C^1$ continuity. When the mass of the wheel is grater than $0$kg (Cases 2 and 4, $\mw = 1000$kg), inertia forces caused by the wheel mass play a role, and, as a result, the second derivative of \Lb\ is required to be twice continuously differentiable (see Equation \eqref{eq:CoupledEq:Direct_2:b}). In these cases, spurious oscillations appear due to the fact that $C^2$ (curvature) continuity is not provided in the bridge model. On the contrary, in Cases 1 and 3 ($\mw = 0$kg) wheels serve merely as geometric constraints between the bridge and the car. As a result, inertia forces associated with the wheel mass do not exist and $C^1$ continuity of the matrix \Lb\ is sufficient (see Equation \eqref{eq:CoupledEq:Direct_2:b}). For Cases 3 and 4, oscillations are more severe due to the additional discontinuity in rotations at the supports of the bridge model. Two approaches are proposed to avoid these oscillations: implementation of a B-spline interpolation (Section~\ref{sec:BSpline}) and an adoption of an alternate time integration scheme (Section~\ref{sec:Bathe}).
\FloatBarrier
\section{Implementation of a B-spline interpolation} \label{sec:BSpline}

In order to avoid spurious oscillations in contact forces and accelerations, a cubic B-spline interpolation is implemented for interpolation of the constraints. Cubic B-spline functions are twice continuously differentiable, which allows to get a $C^2$ interpolating function. Therefore, B-spline functions are used to construct the matrix \Lb.

In this study a B-spline curve is calculated according to \cite{rogers2000introduction}. B-spline is a parametric curve determined by a control polygon.

\begin{equation}
P(t) =\sum _{i=1}^{n+1} B_i N_{i,k}(t) \hskip 1cm t_{\min } \leq t < t_{\max} , \hskip 0.5cm  2 \leq k \leq n+1
\end{equation}

\noindent where $P(t)$ is the position vector along the curve, $B_i$ are the position vectors of the $n+1$ polygon vertices (control points), $N_{i,k}$ are the normalized B-spline basis functions, $t$ is the parameter. The order of B-spline is chosen to be $k=4$. The second derivative of the function $P(t)$ is then continuous, as desired.

The vertices of the control polygon correspond to the nodes of a bridge model.  
Basis functions $N_{i,k}(t)$ are evaluated at each time step for each wheel and depend on the wheels positions along the bridge. The global matrix $\Lb(t)$ is computed using obtained basis functions at each time step. 

\subsection{Numerical example} \label{sec:BSpline:ex}

\begin{figure}[H]   
\centering
\captionsetup{justification=centering}
\captionsetup[subfigure]{justification=centering}
\begin{tabular}{l@{\hskip 0.8cm}r}
\subfloat[Case 2, $\mw = 1000\ kg$. Contact forces (compare with Figure \ref{fig:Bauchau:1:forces:2})] 
{\label{fig:bspline:a}\includegraphics{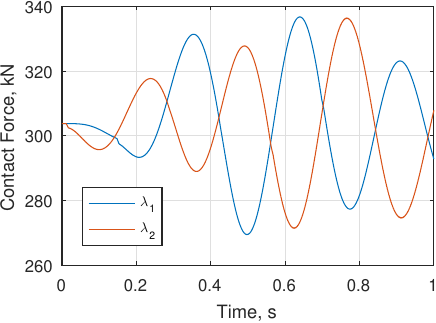}} &
\subfloat[Case 2, $\mw = 1000\ kg$. Acceleration of the midpoint of the first span \hspace{40mm} (compare with Figure~\ref{fig:Bauchau:1:br:2})] 
{\label{fig:bspline:b}\includegraphics{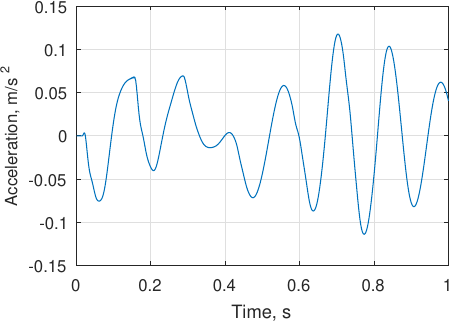}}\\[0pt]
\subfloat[Case 1. Acceleration of the wheels] 
{\label{fig:bspline:c}\includegraphics{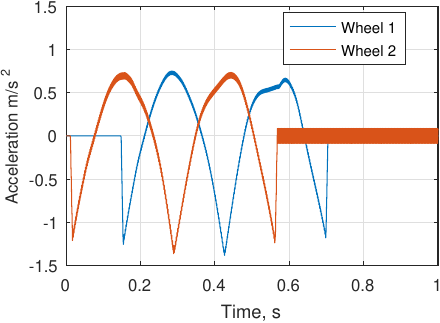}} & 
\subfloat[Case 2. Acceleration of the wheels \hspace{40mm} (compare with Figure~\ref{fig:Bauchau:1:tr:acc:2})] 
{\label{fig:bspline:d}\includegraphics{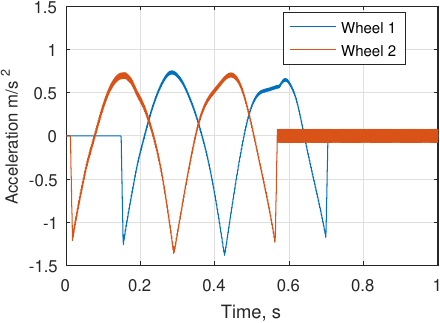}}\\[0pt]
\end{tabular}
\caption{B-spline interpolation. Cases 1 and 2. Contact forces and accelerations of the wheels and the midpoint of the first span of the bridge}
\label{fig:bspline}
\end{figure}

Figures \ref{fig:bspline:b}-\ref{fig:bspline:d} show the obtained accelerations of the wheels and the bridge for Cases 1 and 2. Modification of the matrix \Lb\ considerably reduces spurious oscillations in accelerations, but does not eliminate them completely. Figure \ref{fig:bspline:a} illustrates forces obtained for Case 2. Spurious oscillations in contact forces are not observed. Therefore, implementation of B-spline interpolation in VTSI analysis is reasonable. However, this approach can not be implemented for the bridge system with simply supported ends (Cases 3 and 4) due to the fact that the slope of the spline at the first and last support points should be equal to the slope of the polygon edges. 

\FloatBarrier

\begin{figure}[H]
\centering
\captionsetup{justification=centering}
\includegraphics{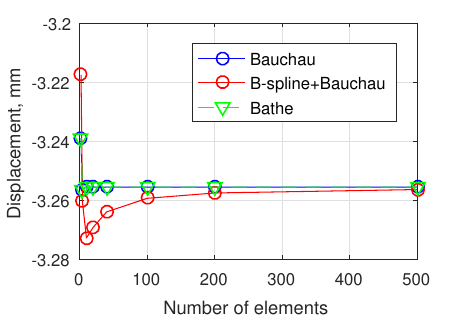}
\caption{Comparison of the convergence of solution for different time integration schemes and different approaches to construction of the matrix \Lb\ }
\label{fig:comparison}
\end{figure}

It should be noted that one of the properties of a B-spline curve is that it does not pass the bridge nodes. As a result, the obtained values of displacements of bridge nodes are approximated and, for a small number of elements $N$, do not correspond to node displacements obtained without using B-spline interpolation. The displacement of the middle node of the first span for Case 1 (see Table \ref{table:cases}) is used to evaluate the convergence of the solution (Figure \ref{fig:comparison}). The figure compares the displacements obtained using the Bauchau integration scheme (Section \ref{sec:Bauchau}), Bauchau integration scheme with B-spline interpolation of the constraints and Bathe integration scheme (Section \ref{sec:Bathe}). As can be seen, when the influence matrix \Lb\ is constructed in a ``standard'' way, that is, by assembling element influence matrices for a current wheel location (Section \ref{sec:2D:bridge}), the solution converges faster. This approach to construction of the matrix \Lb\ was implemented with both, Bauchau and Bathe schemes. However, if the matrix \Lb\ is constructed using B-spline interpolation, as described in current section, the solution converges slower. 

\FloatBarrier

\FloatBarrier
\section{Alternate time integration scheme} \label{sec:Bathe}

Due to the limitations of the B-spline interpolation discussed in Section \ref{sec:BSpline}, an alternate way to mitigate the spurious oscillations is desired. Therefore, in this section the Bathe method is implemented to discretize the coupled equations \eqref{eq:coupled2D} in time.

\subsection{Description of the Bathe method} 

A composite time integration scheme (or Bathe method) was proposed by Bathe~\cite{bathe2005} for solving linear and nonlinear problems in structural dynamics. The main advantage of this scheme is that it remains stable even when the Newmark method produces unstable solutions. The effectiveness of the method in linear analysis was demonstrated by Bathe~\cite{bathe2012}. In this section, the method is summarized for a standard equation of motion, \mbox{${\bf M}{\bf \ddot u} + {\bf C}{\bf \dot u} + {\bf Ku} = {\bf P}$, where ${\bf M}$, ${\bf C}$} and ${\bf K}$ are mass, damping and stiffness matrices, ${\bf P}$ is the vector of applied forces, ${\bf u}$, ${\bf \dot u}$, ${\bf \ddot u}$ are the vectors of nodal displacements (and rotations), velocities and accelerations, respectively. In the next subsection the method is adapted for the case with kinematic constraints.

Assume that the solution is known up to a time point ${n}$. The task is to compute solution at time ${n+1}$. As opposed to majority of time integration schemes, the Bathe method employs calculation of all unknown variables not only at the end (${n+1}$) of a time step \dt, but also at the half (${n+1/2}$) of a time step. For the sake of comparison, in Section~\ref{sec:Bauchau} we obtained the contact forces only at the half (${n+1/2}$) of each time step and all the other quantities at the end (${n+1}$) of each time step. 

According to the Bathe method, in the first sub-step, equations of motions are enforced at time ${n+1/2}$ and the trapezoidal rule (Newmark method with ${\beta =0.25}$ and ${\gamma =0.5}$) is used (equations \eqref{eq:standard:half:b} and \eqref{eq:standard:half:c}). Using notation \eqref{eq:standard:half:subst}, equation \eqref{eq:standard:half:a} can be rewritten as \eqref{eq:standard:half:simpl}.

In the second sub-step, the equations are enforced at time ${n+1}$, and the 3-point backward Euler method is applied (equations \eqref{eq:standard:end:b} and \eqref{eq:standard:end:c}). Using notation \eqref{eq:standard:end:subst}, equation \eqref{eq:standard:end:a} can be rewritten as \eqref{eq:standard:end:simpl}.

\begin{table}[h!]
\footnotesize
\centering
\begin{tabular}{cc} 
\toprule
\textbf{First sub-step}
& \textbf{Second sub-step} \\[0pt]
\midrule 
\parbox[c]{7.8cm}{
    \begin{subequations} \label{eq:standard:half}
    \begin{gather}  
    {\bf M}{{\bf \ddot u}_{n + 1/2}} + {\bf C}{{\bf \dot u}_{n + 1/2}} + {\bf K}{{\bf u}_{n + 1/2}} = {{\bf P}_{n + 1/2}} \label{eq:standard:half:a}\\
    {{\bf \dot u}_{n + 1/2}} = {{\bf \dot u}_n} + \frac{{\Delta t}}{4}\left( {{{\bf \ddot u}_n} + {{\bf \ddot u}_{n + 1/2}}} \right) \label{eq:standard:half:b}\\
    {{\bf u}_{n + 1/2}} = {{\bf u}_n} + \frac{{\Delta t}}{4}\left( {{{\bf \dot u}_n} + {{\bf \dot u}_{n + 1/2}}} \right) \label{eq:standard:half:c}
    \end{gather}  
    \end{subequations}
    } 
& \parbox[c]{7.8cm}{
    \begin{subequations} \label{eq:standard:end}
    \begin{gather}
    {\bf M}{{\bf \ddot u}_{n + 1}} + {\bf C}{{\bf \dot u}_{n + 1}} + {\bf K}{{\bf u}_{n + 1}} = {{\bf P}_{n + 1}} \label{eq:standard:end:a}\\
    {{\bf \dot u}_{n + 1}} = \frac{1}{{\Delta t}}{{\bf u}_n} - \frac{4}{{\Delta t}}{{\bf u}_{n + 1/2}} + \frac{3}{{\Delta t}}{{\bf u}_{n + 1}} \label{eq:standard:end:b}\\
    {{\bf \ddot u}_{n + 1}} = \frac{1}{{\Delta t}}{{\bf \dot u}_n} - \frac{4}{{\Delta t}}{{\bf \dot u}_{n + 1/2}} + \frac{3}{{\Delta t}}{{\bf \dot u}_{n + 1}} \label{eq:standard:end:c}
    \end{gather}  
    \end{subequations}
    }\\[0pt]
\midrule
\parbox[c]{7.8cm}{
    \begin{equation} \label{eq:standard:half:simpl}
    {{\bf \bar M}_{n + 1/2}}{{\bf u}_{n + 1/2}} = {{\bf a}_{n + 1/2}}
    \end{equation}
    }
& \parbox[c]{7.8cm}{
    \begin{equation} \label{eq:standard:end:simpl}
    {{\bf \bar M}_{n + 1}}{{\bf u}_{n + 1}} = {{\bf a}_{n + 1}}
    \end{equation}
    }\\[0pt]
\midrule 
\parbox[c]{7.8cm}{
    \begin{equation} \label{eq:standard:half:subst}
    \begin{gathered}
    {{\bf \bar M}_{n + 1/2}} = \frac{{16}}{{\Delta {t^2}}}{\bf M} + \frac{4}{{\Delta t}}{\bf C} + {\bf K}\\
    {{\bf a}_{n + 1/2}} = {{\bf P}_{n + 1/2}} + {\bf C}\left( {\frac{4}{{\Delta t}}{{\bf u}_n} + {{\bf \dot u}_n}} \right)\\ + {\bf M}\left( {\frac{{16}}{{\Delta {t^2}}}{{\bf u}_n} + \frac{8}{{\Delta t}}{{\bf \dot u}_n} + {{\bf \ddot u}_n}} \right) \\
    \end{gathered}
    \end{equation}
    }
& \parbox[c]{7.8cm}{
    \begin{equation} \label{eq:standard:end:subst}
    \begin{gathered}
    {{\bf \bar M}_{n + 1}} = \frac{9}{{\Delta {t^2}}}{\bf M} + \frac{3}{{\Delta t}}{\bf C} + {\bf K}\\
    {{\bf a}_{n + 1}} = {{\bf P}_{n + 1}} + {\bf C}\left( {\frac{4}{{\Delta t}}{{\bf u}_{n + 1/2}} - \frac{1}{{\Delta t}}{{\bf u}_n}} \right)\\ + {\bf M}\left( \frac{{12}}{{\Delta {t^2}}}{{\bf u}_{n + 1/2}} - \frac{3}{{\Delta {t^2}}}{{\bf u}_n} + \frac{4}{{\Delta t}}{{\bf \dot u}_{n + 1/2}} - \frac{1}{{\Delta t}}{{\bf \dot u}_n} \right)\\
    \end{gathered}
    \end{equation}
    }\\[0pt]
\bottomrule 
\end{tabular}
\caption{Discretization of a standard equation of motion using the Bathe method}
\label{table:Bathe:standard}
\end{table}

\subsection{Adaptation of the Bathe method for VTSI analysis}

Using the Bathe method, equations \eqref{eq:coupled2D} can be discretized in the first sub-step as
\begin{equation} \label{eq:vtsi:half}
\begin{gathered}
\Mt \utddotB_{n + 1/2} + \Ct \utdotB_{n + 1/2} + \Kt \utB_{n + 1/2} + {\Lt}^T \lambdaB_{n + 1/2} = \Pt_{n + 1/2}\\
\Mb \ubddot_{n + 1/2} + \Cb \ubdot_{n + 1/2} + \Kb \ub_{n + 1/2} + \left( \Lb (t_{n + 1/2}) \right)^T \lambdaB_{n + 1/2} = \Pb_{n + 1/2}\\ 
\Lt \utB_{n + 1/2} + \Lb (t_{n + 1/2}) \ub_{n + 1/2} = - \rhoB (t_{n + 1/2})
\end{gathered}
\end{equation}
Using equations \eqref{eq:standard:half:b} and \eqref{eq:standard:half:c} and notation \eqref{eq:standard:half:simpl} and \eqref{eq:standard:half:subst}, equations \eqref{eq:vtsi:half} can be rewritten as \eqref{eq:vtsi:half:2}. Equations \eqref{eq:vtsi:half:2} can then be rewritten again using notation \eqref{eq:vtsi:half:4} to obtain equations~\eqref{eq:vtsi:half:3}.

Similarly to the first sub-step, in the second sub-step we can discretize equations \eqref{eq:coupled2D} as
\begin{equation}  \label{eq:vtsi:end}
\begin{gathered}
\Mt \utddotB_{n + 1} + \Ct \utdotB_{n + 1} + \Kt \utB_{n + 1} + {\Lt}^T \lambdaB_{n + 1} = \Pt_{n + 1}\\
\Mb \ubddot_{n + 1} + \Cb \ubdot_{n + 1} + \Kb \ub_{n + 1} + \left( \Lb (t_{n + 1}) \right)^T \lambdaB_{n + 1} = \Pb_{n + 1}\\
\Lt \utB_{n + 1} + \Lb (t_{n + 1}) \ub_{n + 1} = - \rhoB (t_{n + 1})\\
\end{gathered}
\end{equation}
Using equations \eqref{eq:standard:end:b} and \eqref{eq:standard:end:c} and notation \eqref{eq:standard:end:simpl} and \eqref{eq:standard:end:subst}, equations \eqref{eq:vtsi:end} can be rewritten as \eqref{eq:vtsi:end:2}. Equations \eqref{eq:vtsi:end:2} can then be rewritten again using notation \eqref{eq:vtsi:end:4} to obtain equations~\eqref{eq:vtsi:end:3}.

\begin{table}[h!]
\footnotesize
\centering
\begin{tabular}{cc} 
\toprule
    \textbf{First sub-step}
    & \textbf{Second sub-step} \\[0pt]
\midrule 
    \parbox[c]{7.8cm}{
        \begin{equation} \label{eq:vtsi:half:2}
        \begin{gathered}
        \Mtbar_{n + 1/2} \utB_{n + 1/2} + {\Lt}^T \lambdaB_{n + 1/2} = \at_{n + 1/2}\\
        \Mbbar_{n + 1/2} \ub_{n + 1/2} + (\Lb (t_{n + 1/2}) )^T \lambdaB_{n + 1/2} = \ab_{n + 1/2}\\
        \Lt \utB_{n + 1/2} + \Lb (t_{n + 1/2}) \ub_{n + 1/2} = - \rhoB(t_{n + 1/2})
        \end{gathered}
        \end{equation}
        } 
    & \parbox[c]{7.8cm}{
        \begin{equation} \label{eq:vtsi:end:2}
        \begin{gathered}
        \Mtbar_{n + 1} \utB_{n + 1} + {\Lt}^T \lambdaB_{n + 1} = \at_{n + 1}\\
        \Mbbar_{n + 1} \ub_{n + 1/2} + (\Lb (t_{n + 1}) )^T \lambdaB_{n + 1} = \ab_{n + 1/2}\\
        \Lt \utB_{n + 1} + \Lb (t_{n + 1}) \ub_{n + 1} = - \rhoB(t_{n + 1})
        \end{gathered}
        \end{equation}
        }\\[0pt]
\midrule
    \parbox[c]{7.8cm}{
        \begin{equation} \label{eq:vtsi:half:3}
        \begin{gathered}
        \Mtbar_{n + 1/2} \uttilde_{n + 1/2} = \at_{n + 1/2}\\
        \Mbbar_{n + 1/2} \ubtilde_{n + 1/2} = \ab_{n + 1/2}\\
        \bf{A}_{n + 1/2} \lambdaB_{n + 1/2} = \rhoBbar_{n + 1/2}\\
        \end{gathered}
        \end{equation}
        }
    & \parbox[c]{7.8cm}{
        \begin{equation} \label{eq:vtsi:end:3}
        \begin{gathered}
        \Mtbar_{n + 1} \uttilde_{n + 1} = \at_{n + 1}\\
        \Mbbar_{n + 1} \ubtilde_{n + 1} = \ab_{n + 1}\\
        \bf{A}_{n + 1} \lambdaB_{n + 1} = \rhoBbar_{n + 1}\\
        \end{gathered}
        \end{equation}
        }\\[0pt]
\midrule 
    \parbox[c]{7.8cm}{
        \begin{equation} \label{eq:vtsi:half:4}
        \begin{gathered}
        \uttilde_{n + 1/2} = \utB_{n + 1/2} + \At_{n + 1/2} \lambdaB_{n + 1/2}\\
        \ubtilde_{n + 1/2} = \ub_{n + 1/2} + \Ab_{n + 1/2} \lambdaB_{n + 1/2}\\
        \bf{A}_{n + 1/2} = \left( \Lt \At_{n + 1/2} + \Lb (t_{n + 1/2}) \Ab_{n + 1/2} \right)\\
        \rhoBbar_{n + 1/2} = \Lt \uttilde_{n + 1/2} + \Lb (t_{n + 1/2}) \ubtilde_{n + 1/2} + \rhoB (t_{n + 1/2})\\
        \At_{n + 1/2} = \Mtbar_{n + 1/2} \backslash {\Lt}^T\\
        \Ab_{n + 1/2} = \Mbbar_{n + 1/2} \backslash \left( \Lb (t_{n + 1/2}) \right)^T\\
        \end{gathered}
        \end{equation}
        }
    & \parbox[c]{7.8cm}{
        \begin{equation} \label{eq:vtsi:end:4}
        \begin{gathered}
        \uttilde_{n + 1} = \utB_{n + 1} + \At_{n + 1} \lambdaB_{n + 1}\\
        \ubtilde_{n + 1} = \ub_{n + 1} + \Ab_{n + 1} \lambdaB_{n + 1}\\
        \bf{A}_{n + 1} = \left( \Lt \At_{n + 1} + \Lb (t_{n + 1}) \Ab_{n + 1} \right)\\
        \rhoBbar_{n + 1} = \Lt \uttilde_{n + 1} + \Lb (t_{n + 1}) \ubtilde_{n + 1} + \rhoB (t_{n + 1})\\
        \At_{n + 1} = \Mtbar_{n + 1} \backslash {\Lt}^T\\
        \Ab_{n + 1} = \Mbbar_{n + 1} \backslash \left( \Lb (t_{n + 1}) \right)^T\\
        \end{gathered}
        \end{equation}
        }\\[0pt]
\bottomrule 
\end{tabular}
\caption{Discretization of equations of motion \eqref{eq:coupled2D} using the Bathe method}
\label{table:Bathe:vtsi}
\end{table}

\subsection{Modified procedure} 

Procedure~\ref{procedure_2D} is modified in order to replace the Bauchau time integration scheme (Section~\ref{sec:Bauchau}) with the Bathe method. The updated procedure is presented below. 

\floatname{algorithm}{Procedure} 
\begin{algorithm}[H] 
\caption{\text Linear dynamic analysis of bridge-train interaction using the Bathe method} \label{procedure_2D_Bathe}
\begin{algorithmic}[1]
\State Given: \Pt, \Pb, $\rhoB(t)$, \dt 
\State Assemble \Mt, \Ct, \Kt
\State Assemble $\Mtbar_{n+1/2}$ and $\Mbbar_{n+1/2}$ 
\State Assemble \Mb, \Cb, \Kb
\State Assemble $\Mbbar_{n+1/2}$ and $\Mtbar_{n+1}$ 
\State Compute $\utB_0$ by linear static analysis. Set $\utdotB_0 = {\bf 0}$ \Comment{Train initial conditions}
\State Compute $\ub_0$ by linear static analysis. Set $\ubdot_0 = {\bf 0}$ \Comment{Bridge initial conditions}
\For{$n\leftarrow 0$, NUMINC} 
\State Compute $\xw(t_{n+1/2})$ and $\xw(t_{n+1})$ by equation \eqref{wheelpos}
\State Obtain $\Lb(t_{n+1/2})$ and $\Lb(t_{n+1})$

\State \textbf{\emph{First sub-step}} \hrulefill
\State Compute $\at_{n+1/2}$ and $\ab_{n+1/2}$ 
\State Solve $\Mtbar_{n+1/2} \left[ \uttilde_{n+1/2}\; \At_{n+1/2} \right] = \left[ \at_{n+1/2}\; {\Lt}^T \right]$
\State Solve $\Mbbar_{n+1/2} \left[ \ubtilde_{n+1/2}\; \Ab_{n+1/2} \right] = \left[ \ab_{n+1/2}\; \Lb \left( t_{n + 1/2} \right)^T \right]$
\State Solve ${\bf A}_{n + 1/2} \lambdaB_{n + 1/2} = \rhoBbar_{n + 1/2}$
\State Set $\utB_{n+1/2} = \uttilde_{n+1/2} - \At_{n+1/2} \lambdaB_{n + 1/2}$; $\utdotB_{n + 1/2} = \left( \utB_{n + 1/2} - \utB_n \right) \frac{4}{\dt} - \utdotB_n$; \hfill\break \hspace*{1.5cm} $\utddotB_{n + 1/2} = \left( \utdotB_{n + 1/2} - \utdotB_n \right) \frac{4}{\dt} - \utddotB_n$;
\State Set $\ub_{n+1/2} = \ubtilde_{n+1/2} - \Ab_{n+1/2} \lambdaB_{n + 1/2}$;  $\ubdot_{n + 1/2} = \left( \ub_{n + 1/2} - \ub_n \right) \frac{4}{\dt} - \ubdot_n$;  \hfill\break \hspace*{1.5cm} $\ubddot_{n + 1/2} = \left( \ubdot_{n + 1/2} - \ubdot_n \right) \frac{4}{\dt} - \ubddot_n$

\State \textbf{\emph{Second sub-step}} \hrulefill
\State Compute $\at_{n+1}$ and $\ab_{n+1}$  
\State Solve $\Mtbar_{n+1} \left[ \uttilde_{n+1}\; \At_{n+1} \right] = \left[ \at_{n+1}\; {\Lt}^T \right]$
\State Solve $\Mbbar_{n+1} \left[ \ubtilde_{n+1}\; \Ab_{n+1} \right] = \left[ \ab_{n+1}\; \Lb \left( t_{n + 1} \right)^T \right]$
\State Solve ${\bf A}_{n + 1} \lambdaB_{n + 1} = \rhoBbar_{n + 1}$
\State Set $\utB_{n+1} = \uttilde_{n+1} - \At_{n+1} \lambdaB_{n + 1}$;        $\utdotB_{n + 1} = \frac{1}{\dt}\utB_n - \frac{4}{\dt}\utB_{n + 1/2} + \frac{3}{\dt}\utB_{n + 1}$; \hfill\break \hspace*{1.5cm} $\utddotB_{n + 1} = \frac{1}{\dt} \utdotB_n - \frac{4}{\dt}\utdotB_{n + 1/2} + \frac{3}{\dt} \utdotB_{n + 1}$;
\State Set $\ub_{n+1} = \ubtilde_{n+1} - \Ab_{n+1} \lambdaB_{n + 1}$;         $\ubdot_{n + 1} = \frac{1}{\dt} \ub_n - \frac{4}{\dt} \ub_{n + 1/2} + \frac{3}{\dt} \ub_{n + 1}$; \hfill\break \hspace*{1.5cm} $\ubddot_{n + 1} = \frac{1}{\dt} \ubdot_n - \frac{4}{\dt} \ubdot_{n + 1/2} + \frac{3}{\dt} \ubdot_{n + 1}$
\EndFor
\end{algorithmic}
\end{algorithm}

\FloatBarrier
\subsection{Numerical example} \label{sec:Bathe:ex}

\begin{figure}[H]   
\centering
\captionsetup{justification=centering}
\captionsetup[subfigure]{justification=centering}
\begin{tabular}{l@{\hskip 0.8cm}r}
\subfloat[Case 2. Contact forces] 
{\label{fig:Bathe:1:a}\includegraphics{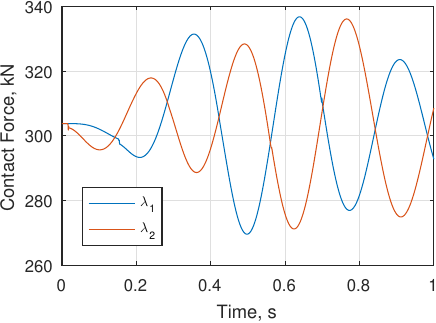}} &
\subfloat[Case 2. Acceleration of the midpoint of the first span] 
{\label{fig:Bathe:1:b}\includegraphics{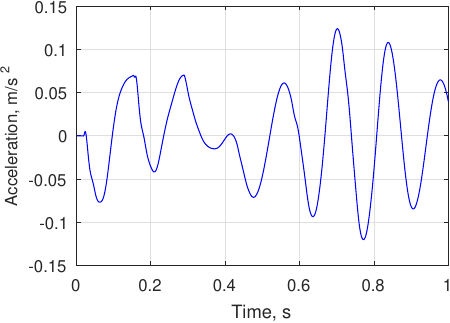}}\\
\subfloat[Case 1. Acceleration of the wheels] 
{\label{fig:Bathe:1:c}\includegraphics{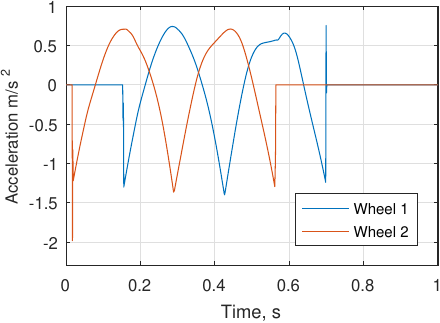}} & 
\subfloat[Case 2. Acceleration of the wheels] 
{\label{fig:Bathe:1:d}\includegraphics{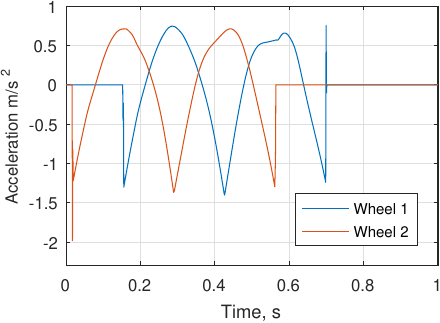}}\\[0pt]
\end{tabular}
\caption{Bathe method. Cases 1 and 2. Contact forces and accelerations of the wheels and the midpoint of the first span of the bridge}
\label{fig:Bathe:1}
\end{figure}

\begin{figure}[H]   
\centering
\captionsetup{justification=centering}
\captionsetup[subfigure]{justification=centering}
\begin{tabular}{l@{\hskip 0.8cm}r}
\subfloat[Case 4. Contact forces] 
{\label{fig:Bathe:0:a}\includegraphics{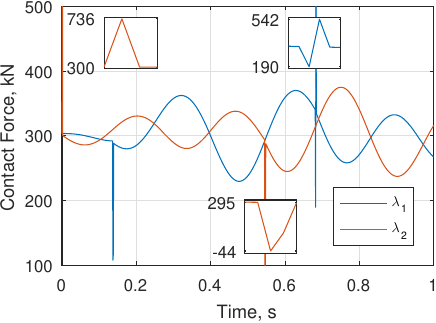}} &
\subfloat[Case 3. Acceleration of the midpoint of the first span] 
{\label{fig:Bathe:0:b}\includegraphics{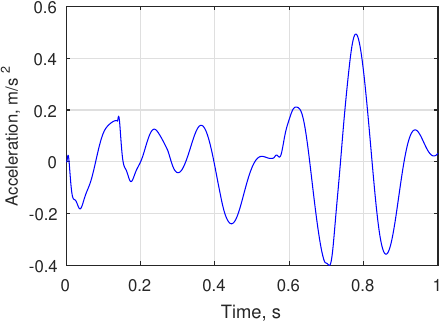}}\\[0pt]
\subfloat[Case 3. Acceleration of the wheels] 
{\label{fig:Bathe:0:c}\includegraphics{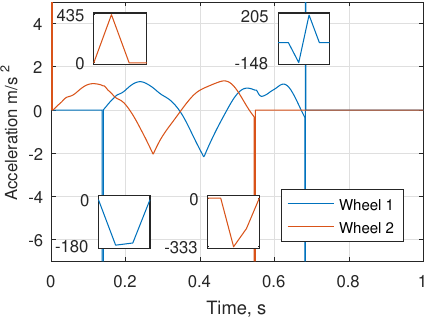}} & 
\subfloat[Case 4. Acceleration of the wheels] 
{\label{fig:Bathe:0:d}\includegraphics{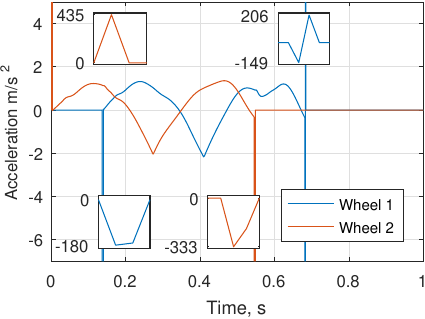}}\\[0pt]
\end{tabular}
\caption{Bathe method. Cases 3 and 4. Contact forces and accelerations of the wheels and the midpoint of the first span of the bridge}
\label{fig:Bathe:0}
\end{figure}

In this section Cases 1 to 4 (see Table \ref{table:cases}) are analyzed using the Bathe method. Figures \ref{fig:Bathe:1} and \ref{fig:Bathe:0} present the obtained results.

As can be seen, the implementation of the Bathe method mitigates the spurious oscillations in the accelerations and contact forces. However, spikes are observed in the wheels accelerations (Figures \ref{fig:Bathe:1:c}, \ref{fig:Bathe:1:d}, and specially \ref{fig:Bathe:0:c}, \ref{fig:Bathe:0:d}) and the contact forces (Figure \ref{fig:Bathe:0:a}). These spikes occur in the time steps when the wheels are entering or leaving the bridge due to the fact that curvature and rotation discontinuity trigger high frequency oscillations (as discussed in Section \ref{sec:Bauchau:ex}). Similar ``undershoots'' were observed in the results of an example analysis performed by Bathe \cite{bathe2012}. But the Bathe integration scheme attenuates these oscillations very quickly (with one time step). Therefore, no further special measures are taken here regarding these spikes. If so desired, the transition can be smoothed just at the start and end of the bridge.

The stability of the obtained solutions can be explained by the fact that the Bathe method includes some numerical damping. This is due to the employment of the backward Euler method \eqref{eq:standard:end} in the second sub-step \cite{bathe2005}.

\FloatBarrier
\subsubsection{Verification of the algorithm using SundialsTB} \label{sec:Sundials}

In order to verify the obtained algorithm, we solve the system of equations \eqref{eq:coupled2D} using IDAS, a solver for systems of differential-algebraic equations (DAE). IDAS (Implicit Differential-Algebraic solver with
Sensitivity capabilities) is a part of a software suite SUNDIALS (SUite of Nonlinear and DIfferential/ALgebraic
equation Solvers) \cite{serban2015user}. In order to utilize earlier developed MATLAB functions for bridge, train and constraints modeling, we employ SundialsTB \cite{SUNDIALS}, a MATLAB interface to the SUNDIALS.

The system \eqref{eq:coupled2D} is a DAE system of the form

\begin{equation}
\begin{gathered}
{\bf F} \left( t, {\bf y}, {\bf \dot y} \right) = 0,\\
{\bf y}(t_0) = {\bf y}_0, {\bf \dot y}(t_0) = {\bf \dot y}_0
\end{gathered}
\end{equation} 

\noindent where ${\bf y} = {\left[ {\bf u}^t\ {\bf u}^b\ {\bf v}^t\ {\bf v}^b \right] }^T$, containing vectors of train displacements and velocities, bridge displacements and velocities respectively. We set initial values of displacements based on the separate static analysis for the bridge and the train. Initial values of velocities are assigned to be zero. IDAS implements BDF (Backward
Differentiation Formula) integration method. 
Residual functions ${\bf F}$ are computed as

\begin{equation} \label{residuals}
\begin{aligned}
{{\bf F}_1} & = {{\bf M}^t}{{\bf \dot v}^t} + {{\bf C}^t}{{\bf v}^t} + {{\bf K}^t}{{\bf u}^t} + {\left( {{{\bf L}^t}} \right)^T}{\boldsymbol \lambda} - {{\bf P}^t} \\
{{\bf F}_2} & = {{\bf M}^b}{{\bf \dot v}^b} + {{\bf C}^b}{{\bf v}^b} + {{\bf K}^b}{{\bf u}^b} + {\left( {{{\bf L}^b}(t)} \right)^T}{\boldsymbol \lambda} - {{\bf P}^b} \\
{{\bf F}_3} & = {{\bf v}^t} - {{\bf \dot u}^t} \\
{{\bf F}_4} & = {{\bf v}^b} - {{\bf \dot u}^b} \\
{{\bf F}_5} & = {{\bf L}^t}{{\bf u}^t} + {{\bf L}^b}(t){{\bf u}^b} - {\boldsymbol \rho} \left( t \right) 
\end{aligned}
\end{equation}

\begin{figure}[H]  
\centering
\captionsetup{justification=centering}
\captionsetup[subfigure]{justification=centering}
\begin{tabular}{l@{\hskip 0.8cm}r}
\subfloat[Contact forces] 
{\label{fig:SUN:a}\includegraphics{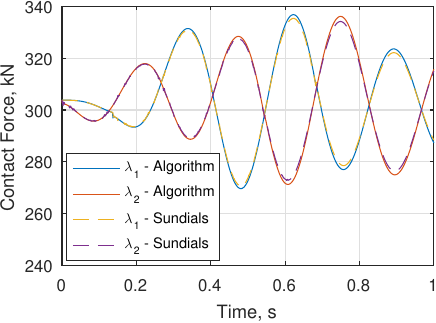}} &
\subfloat[Displacement and rotation of the middle point of the first span] 
{\label{fig:SUN:b}\includegraphics{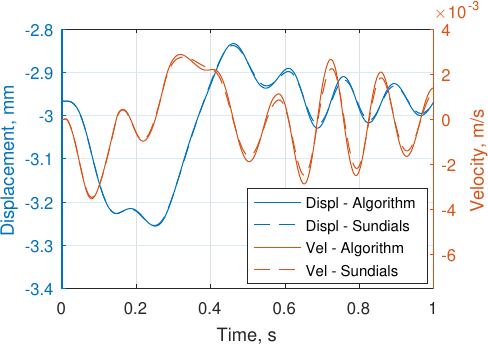}}\\[0pt]
\end{tabular}
\caption{Comparison of the solutions obtained by Bathe method and using SundialsTB. Contact forces, displacement and rotation of the middle point of the first span of the bridge}
\label{fig:SUN}
\end{figure}

\begin{figure}[H]  
\centering
\captionsetup{justification=centering}
\captionsetup[subfigure]{justification=centering}
\begin{tabular}{l@{\hskip 0.8cm}r}
\subfloat[Displacements of the wheels] 
{\label{fig:SUN:tr:a}\includegraphics{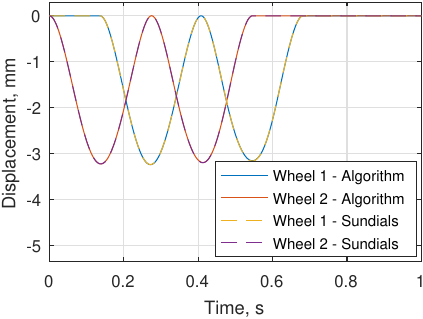}} &
\subfloat[Displacements and rotations of the car] 
{\label{fig:SUN:tr:b}\includegraphics{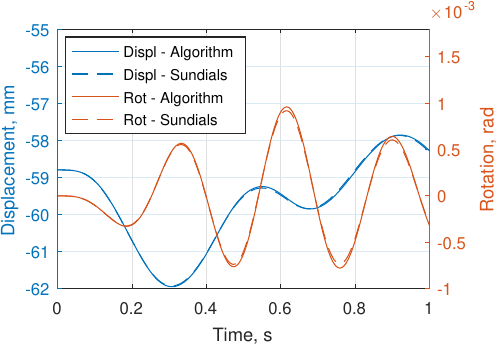}}\\[0pt]
\subfloat[Velocities of the wheels] 
{\label{fig:SUN:tr:c}\includegraphics{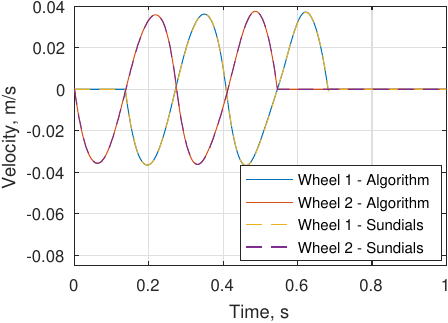}} &
\subfloat[Velocities and rotational velocities of the car] 
{\label{fig:SUN:tr:d}\includegraphics{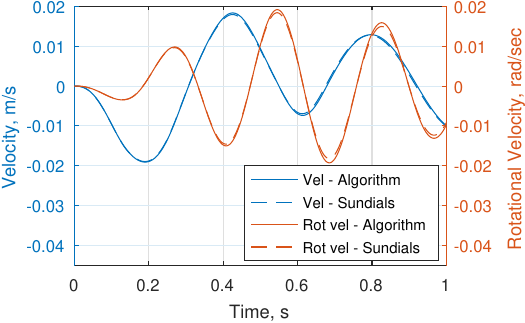}}\\[0pt]
\end{tabular}
\caption{Comparison of the solutions obtained by Bathe method and using SundialsTB. Displacements and velocities of the train model}
\label{fig:SUN:tr}
\end{figure}

A solution for Case 1 (see Table \ref{table:cases}) is computed using SundialsTB. Figures \ref{fig:SUN} and \ref{fig:SUN:tr} illustrate a comparison of the results obtained using the Bathe method with the SundialsTB solution of the system \eqref{residuals}. The results are in good agreement, which means that system \eqref{eq:coupled2D} was discretized correctly using the Bathe method. Therefore, this time integration scheme is chosen to be used for VTSI algorithm.

\FloatBarrier
\subsubsection{Modeling of track irregularities} \label{sec:irreg}

The following power spectral density (PSD) function is used to represent the elevation (vertical) track irregularities \cite{Claus1998}:
\begin{equation}
   S(\Omega_n) = \dfrac{A \Omega_c^2}{(\Omega_n^2 + \Omega_r^2)(\Omega_n^2 + \Omega_c^2)}
\end{equation}
where $\Omega_n = n \Delta \Omega = n (\Omega_u - \Omega_l) /N,\ n = 1,...,N-1$ are discrete frequencies, $\Omega_r = 0.0206 {\text{rad/m}}$ and $\Omega_c = 0.825 {\text{rad/m}}$ are coefficients. The value of coefficient $A$ is set to $1.5 \cdot 10^{-6}$, which corresponds to track Class 6 according to the Federal Railroad Administration (FRA) classification \cite{yang2004vehicle}. The uppermost frequency is set to $\Omega_u = 13.57383 {\text{rad/m}}$, the lowest frequency is equal to $\Omega_l = 0.00383 {\text{rad/m}}$  and the total number of frequencies is chosen to be $N = 3540$ \cite{Claus1998}.

\begin{figure}[H]  
\centering
\captionsetup{justification=centering}
\captionsetup[subfigure]{justification=centering}
\begin{tabular}{c l}
\subfloat[Numerical example: Two-span bridge] 
    {\hskip 0.76cm \label{fig:irreg:a}\includegraphics[scale=0.27]{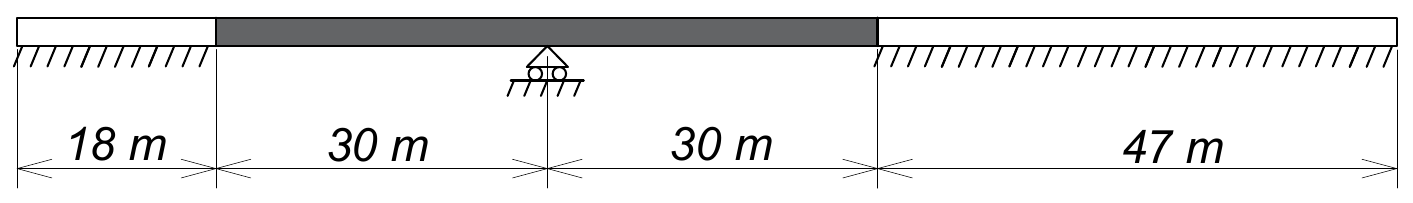}} 
& \\[0pt]
\subfloat[Rail irregularity profile corresponding to Figure \ref{fig:irreg:a}] 
    {\label{fig:irreg:b}\includegraphics{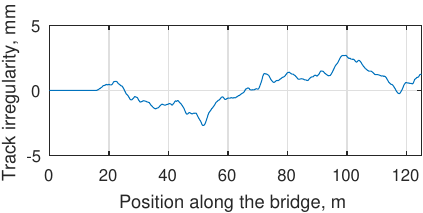}} 
& \multirow{2}{*}[13em]{ \subfloat[Contact forces] 
    {\label{fig:irreg:c}\includegraphics{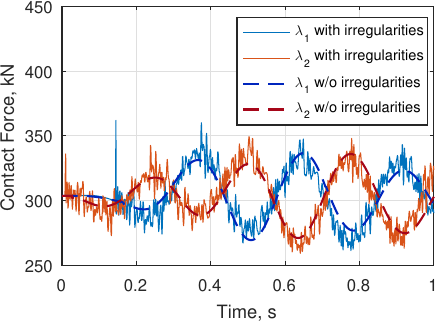}}}\\
\subfloat[Acceleration of the middle point of the first span] 
    {\label{fig:irreg:d}\includegraphics{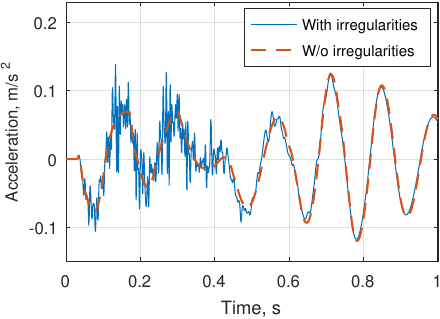}}
& \subfloat[Acceleration of the car] 
    {\label{fig:irreg:e}\includegraphics{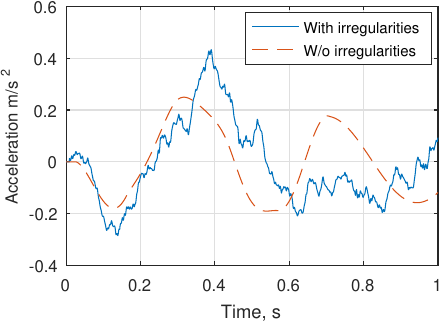}}\\[0pt]
\end{tabular}
\caption{Case 2. (a) The analysed bridge, (b) corresponding track irregularity profile, (c) contact forces between the train and the bridge, responses of (d) the bridge and (e) the car}
\label{fig:irreg}
\end{figure}

Using the spectral representation method, the track irregularities can be generated as a function of a position $x$ along the bridge
\begin{equation}
    \rho(x) = \sqrt{2} \sum_{n = 1}^{N-1} A_n cos(\Omega_n x + \varphi_n)
\end{equation}
where $\varphi_n,\ n = 1,...,N-1$ are random phase angles uniformly distributed in the range $[0,2\pi]$. The coefficients $A_n$ can be found as $A_n = \sqrt{\del{ \dfrac{1}{\pi} S(\Omega_n) + a_n S(0)} \Delta \Omega}$, where $a_1 = \dfrac{4}{6\pi} $, $a_2 = \dfrac{1}{6\pi}$ and $a_n = 0$ for $n = 3,...,N-1$. 

The generated track irregularity profile is normalized so that deviations do not exceed the maximum tolerable deviations of $2.7$mm for track Class 6 \cite{yang2004vehicle}. The normalized irregularity profile is presented in Figure \ref{fig:irreg:b} and applied to two-span continuous bridge (Figures \ref{fig:irreg:a} and \ref{fig:cases:1-4}). It should be noted that track irregularities were set to zero at the beginning of the track (outside of the first bridge span) in order to satisfy initial conditions. Case 2 (see Table \ref{table:cases}) is then analysed. As can be seen from Figures \ref{fig:irreg:c} - \ref{fig:irreg:e}, responses of the bridge and the car are amplified if track irregularities are presented. 

\FloatBarrier
\subsubsection{Influence of train speed on bridge response}\label{sec:resonance}

Case 5 (see Table \ref{table:cases} and Figure \ref{fig:cases:5}) is analyzed in this section in order to demonstrate the influence of the train speed on the response of the bridge. The properties of an one-span simply supported bridge correspond to the properties of bridges in Cases 1-4 (see Section~\ref{sec:Bauchau:ex}). The parameters of a single train car correspond to the car parameters from Cases 1-4 as well; the train is composed of 10 identical cars.

\begin{figure}[H]  
\centering
\captionsetup{justification=centering}
\captionsetup[subfigure]{justification=centering}
\begin{tabular}{l r}
\subfloat[Maximum absolute displacement of the midspan with respect to the train speed] 
{\label{fig:res:a}\includegraphics{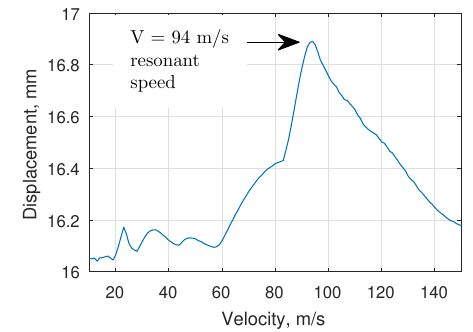}}
& \subfloat[Maximum absolute acceleration of the midspan with respect to the train speed] 
{\label{fig:res:b}\includegraphics{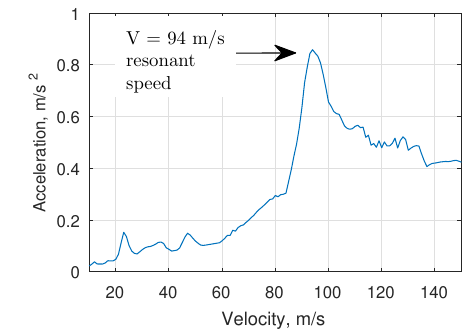}}\\[0pt]
\subfloat[Displacement of the midspan corresponding to the train speed $94 {\text{m/s}}$] 
{\hskip 0.6cm \label{fig:res:c}\includegraphics{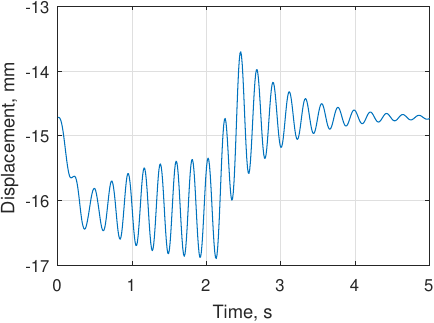}}
& \subfloat[Acceleration of the midspancorresponding to the train speed $94 {\text{m/s}}$] 
{\label{fig:res:d}\includegraphics{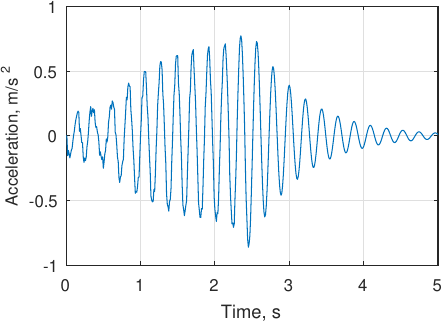}}\\[0pt]
\end{tabular}
\caption{Case 5. Response of the midspan of the bridge: (a) and (b) depending on the train speed; (c)~and~(d)~corresponding to the resonant train speed of $94 {\text{m/s}}$}
\label{fig:res}
\end{figure}

According to Yang~et.~al~\cite{yang2004vehicle}, the resonant train speed can be found as $v^{\text{t}}_{\text{res}} = f^{\text{b}} \lct$, where $f^{\text b}$ is the natural frequency of a~bridge and $\lct$ is the vehicle length. 
The first natural frequency of the bridge is $f^{\text b} = 4.6$Hz. The resonant speed can then be predicted to be $v^{\text{t}}_{\text{res}} = 92{\text{m/s}}$, which is close to $94 {\text{m/s}}$, the~resonant speed obtained through the numerical analysis of Case 5 for the train speed range $[10,150]{\text{m/s}}$ (Figures \ref{fig:res:a} and \ref{fig:res:b}). The effect of resonance is illustrated in Figures \ref{fig:res:c} and \ref{fig:res:d}, which show a significant amplification of the bridge response as 10 cars passing the bridge with the resonant speed of $94{\text{m/s}}$.

\FloatBarrier
\section{Modeling of contact separation} \label{sec:separation}

\begin{figure}[h]
\centering
\captionsetup{justification=centering}
\begin{tabular}{l@{\hskip 0.5cm}r}
\subfloat[Contact separation between the wheel and the track] 
{\label{fig:constr_inequal:real} {\includegraphics[scale=0.27]{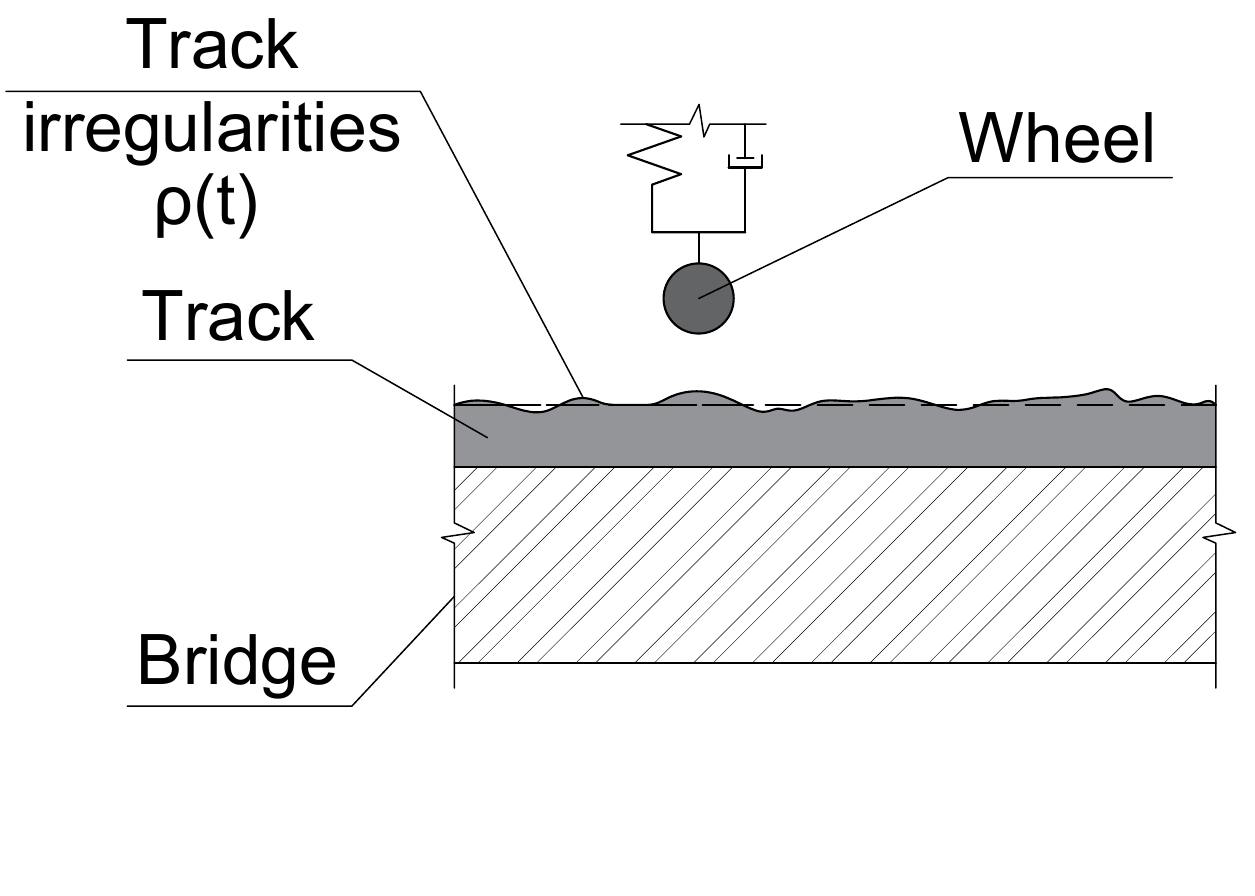}}} &
\subfloat[Simplified modeling of contact separation] 
{\label{fig:constr_inequal:model} \includegraphics[scale=0.3]{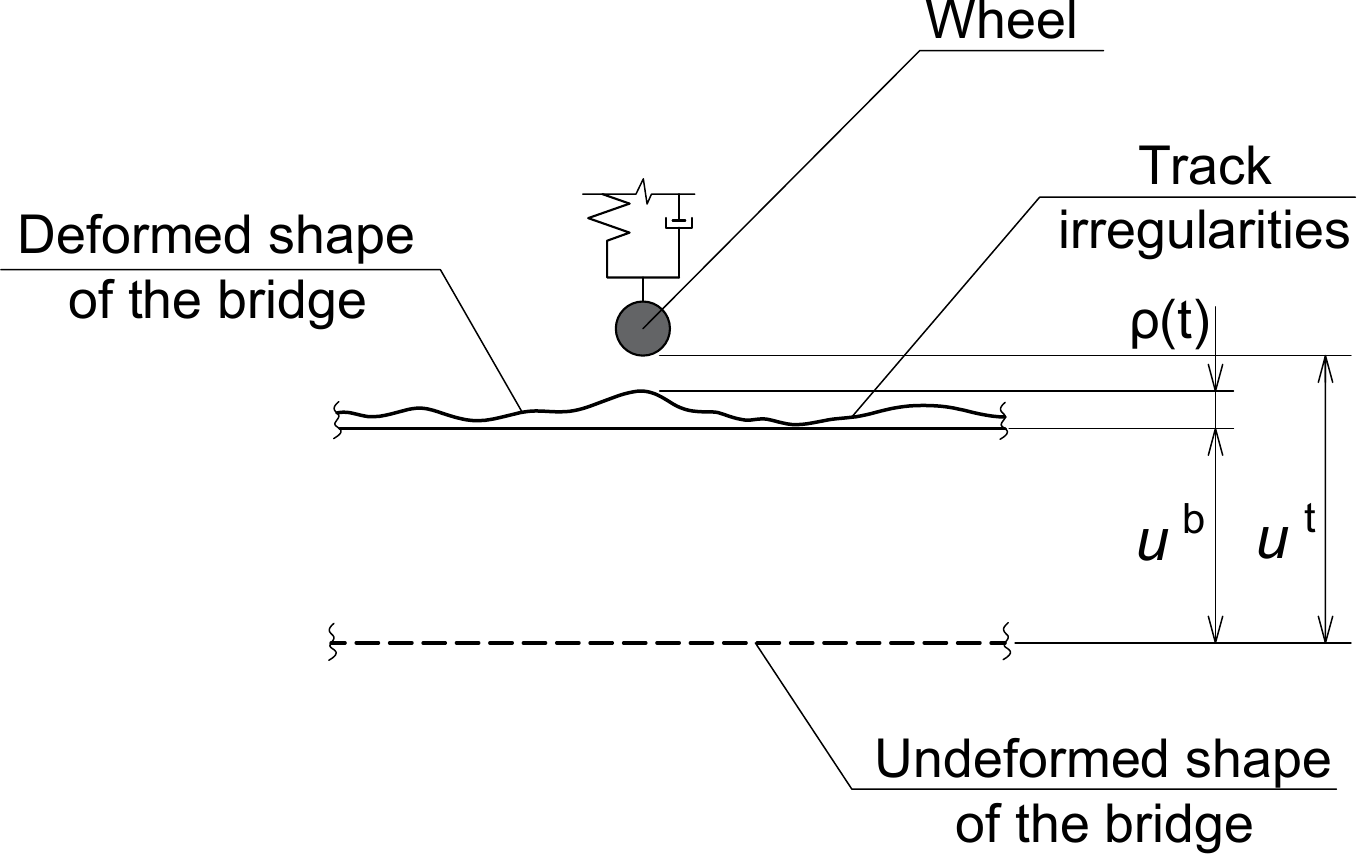}}
\end{tabular}
\caption{Modeling of contact separation}
\label{fig:constr_inequal}
\end{figure}

To model the separation between the wheels and the bridge (Figure \ref{fig:constr_inequal}), the equality constraint \eqref{eq:coupled2Dc} is replaced with an inequality constraint. Combining the obtained inequality with equations of motion \eqref{eq:coupled2Da} and \eqref{eq:coupled2Db}, we obtain the coupled system of equations

\begin{subequations} \label{sep_coupled2D}
\begin{gather}  
\Mt \utddotB + \Ct \utdotB + \Kt \utB + \left( \Lt \right)^T \lambdaB = \Pt \label{sep_coupled2D_a}\\
\Mb \ubddot + \Cb \ubdot + \Kb \ub + \left( \Lb(t) \right)^T \lambdaB = \Pb \label{sep_coupled2D_b}\\
\Lt \utB + \Lb(t) \ub \leq - \rhoB \left( t \right) \label{sep_coupled2D_c}
\end{gather} 
\end{subequations}

\FloatBarrier
\FloatBarrier

Then, we have the following linear complementarity conditions \cite{cottle1992linear}

\begin{subequations} \label{sep_constr}
\begin{gather}  
\Lt \utB + \Lb (t) \ub \le - \rhoB (t) \label{sep_constr_a}\\
\lambdaB  \ge 0 \label{sep_constr_b}\\
\lambdaB^T \left( \Lt \utB + \Lb (t) \ub + \rhoB (t) \right) = {\bf 0} \label{sep_constr_c}
\end{gather} 
\end{subequations}

Equation \eqref{sep_constr_a} implies that the wheel may stay on the bridge or separate from the bridge, but penetrations are not possible. According to equation \eqref{sep_constr_b}, contact force \lambdaB\ is always non-negative: it is positive when the wheel is on the bridge and wheel force is applied on the bridge, and equal to zero when the wheel is separated from the bridge. Equation \eqref{sep_constr_c} implies that at least one one the components of the equation must be equal to zero at each time step, that is, the contact force is zero when separation occurs, and when contact force is positive, there is no separation gap between the wheel and the bridge.

Equation \eqref{sep_coupled2D_c} is discretized in time using the Bathe method. Using equations \eqref{eq:vtsi:half} and \eqref{eq:vtsi:end}, discretized equations in the first and the second sub-step are written as

\begin{equation} \label{sep_discr}
\begin{gathered}
\Lt \utB_{n + 1/2} + \Lb (t_{n + 1/2}) \ub_{n + 1/2} \le - \rhoB (t_{n + 1/2})\\
\Lt \utB_{n + 1} + \Lb \left( t_{n + 1} \right) \ub_{n + 1} \le  - \rhoB \left( t_{n + 1} \right)\\
\end{gathered}
\end{equation}

Using equations \eqref{eq:vtsi:half} and \eqref{eq:vtsi:end}, equations \eqref{sep_discr} can be rewritten as
\begin{equation} \label{sep_discr_2}
\begin{gathered}
\left( \Lt \At_{n + 1/2} + \Lb(t_{n + 1/2})\Ab_{n + 1/2} \right) \lambdaB_{n + 1/2} \ge \Lt \uttilde_{n + 1/2} + \Lb(t_{n + 1/2})\ubtilde_{n + 1/2} + \rhoB (t_{n + 1/2})\\
\left( \Lt \At_{n + 1} + \Lb(t_{n + 1})\Ab_{n + 1} \right) \lambdaB_{n + 1} \ge \Lt \uttilde_{n + 1} + \Lb(t_{n + 1}) \ubtilde_{n + 1} + \rhoB (t_{n + 1})
\end{gathered}
\end{equation}
or concisely, using notation of Procedure \ref{procedure_2D_Bathe}, as
\begin{equation} \label{sep_discr_3}
\begin{gathered}
{\bf A}_{n + 1/2} \lambdaB_{n + 1/2} \ge \rhoBbar_{n + 1/2}\\
{\bf A}_{n + 1} \lambdaB_{n + 1} \ge \rhoBbar_{n + 1}
\end{gathered}
\end{equation}

Then, considering inequalities \eqref{sep_discr_3} and constraints \eqref{sep_constr}, we can formulate the contact separation problem as a linear complementarity problem (LCP) for both sub-steps
\begin{subequations} \label{LCP_discr}
\begin{gather} 
\begin{gathered}
{\bf A}_{n + 1/2} \lambdaB_{n + 1/2} - \rhoBbar_{n + 1/2} \ge {\bf 0},\ \lambdaB_{n + 1/2} \ge {\bf 0},\\ \lambdaB_{n + 1/2}^T \left( \Lt \utB_{n+1/2} + \Lb(t_{n+1/2})\ub_{n+1/2} + \rhoB (t_{n+1/2}) \right) = {\bf 0} \label{LCP_discr_a}\\
\end{gathered}
\end{gather}
\begin{gather} 
\begin{gathered}
{\bf A}_{n + 1} \lambdaB_{n + 1} - \rhoBbar_{n + 1} \ge {\bf 0},\ \lambdaB_{n + 1} \ge {\bf 0},\\ \lambdaB_{n + 1}^T \left( \Lt \utB_{n+1} + \Lb (t_{n+1}) \ub_{n+1} + \rhoB (t_{n+1}) \right) = {\bf 0} \label{LCP_discr_b}
\end{gathered}
\end{gather} 
\end{subequations}
The goal is to find the vector $\lambdaB_{n+1/2}$ satisfying constraints \eqref{LCP_discr_a} in the first sub-step and the vector $\lambdaB_{n+1}$ satisfying constraints \eqref{LCP_discr_b} in the second sub-step.

\subsection{Modified procedure}

To solve the LCP \eqref{LCP_discr} we use a pivot algorithm developed by Almqvist et al. \cite{LCP} for solving LCPs.

Procedure \ref{procedure_2D_Bathe} is modified to include a call to the LCP solver into the code. The modified part of the Procedure \ref{procedure_2D_Bathe} is presented below. The line numbers in Procedure \ref{procedure_2D_LCP} correspond to line numbers in Procedure \ref{procedure_2D_Bathe}. As can be noted, only two lines of the Procedure \ref{procedure_2D_Bathe} need to be changed in order to enable the possibility of contact separation between the wheels and the bridge.

\floatname{algorithm}{Procedure} 
\begin{algorithm} 
\caption{\text Linear dynamic analysis of bridge-train interaction using the Bathe method and enabling contact separation} \label{procedure_2D_LCP}
\begin{algorithmic}[1]
\Statex . . .

\makeatletter                 
\setcounter{ALG@line}{7}      
\makeatother                  
\For{$n\leftarrow 0$, NUMINC}
\Statex \hspace*{0.5cm} . . .

\makeatletter                 
\setcounter{ALG@line}{10}     
\makeatother                  
\State \textbf{\emph{First sub-step}} \hrulefill
\Statex \hspace*{0.5cm} . . .

\makeatletter                 
\setcounter{ALG@line}{16}     
\makeatother                  
\State Solve $\lambdaB_{n + 1/2} = \text{LCPSolve} \left( {\bf A}_{n + 1/2},- \rhoBbar_{n + 1/2} \right) $
\Statex \hspace*{0.5cm} . . .

\makeatletter                 
\setcounter{ALG@line}{19}     
\makeatother                  
\State \textbf{\emph{Second sub-step}} \hrulefill
\Statex \hspace*{0.5cm} . . .

\makeatletter                 
\setcounter{ALG@line}{25}     
\makeatother                  
\State Solve $\lambdaB_{n + 1} = \text{LCPSolve} \left( {\bf A}_{n + 1},- \rhoBbar_{n + 1} \right) $
\Statex \hspace*{0.5cm} . . .

\EndFor
\end{algorithmic}
\end{algorithm}

\FloatBarrier
\subsection{Numerical example} \label{sec:LCP:ex}

In this section, Case 6 from Table \ref{table:cases} is analyzed to investigate the performance of the obtained scheme with contact separation enabled.

The bridge is represented by two-span continuous beam. The ends on the bridge are fixed.\\
The material and geometric properties of the bridge are as follows: length of one span $25$m, Young's modulus ${E=22}GPa$, cross-section second moment of inertia ${I_y=4 {\text m}^4}$, mass per unit length ${38000 {\text{kg/m}}}$. The damping is modeled using Rayleigh model as described in Section \ref{sec:Bauchau}.

The train consists of two cars with masses ${m^c=100000}$kg, wheel-to-wheel length ${l^c=6}$m, moment of inertia ${I^c=506670 \text{kg} \cdot \text{m}^2}$. The distance between the last wheel of the first car and the first wheel of the second car is equal to $3$m. Each suspension has stiffness 

\begin{figure}[H] 
\centering
\captionsetup{justification=centering}
\captionsetup[subfigure]{justification=centering}
\begin{tabular}{l@{\hskip 1.2cm}r}
\subfloat[] 
{\label{fig:LCP:forces}\includegraphics{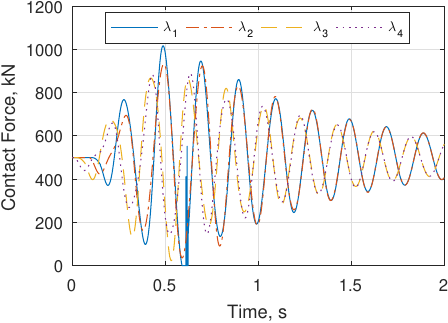}} &
\subfloat[] 
{\label{fig:LCP:displ}\includegraphics{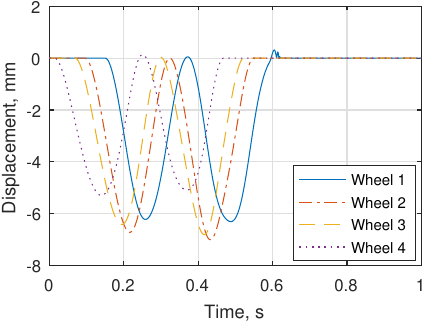}}
\end{tabular}
\caption{Case 6. Contact forces (a) and displacements of the wheels (b)}
\label{fig:LCP}
\end{figure}

\begin{figure}[H] 
\centering
\captionsetup{justification=centering}
\captionsetup[subfigure]{justification=centering}
\begin{tabular}{l@{\hskip 1.2cm}r}
\subfloat[] 
{\label{fig:LCP:detail:forces}\includegraphics{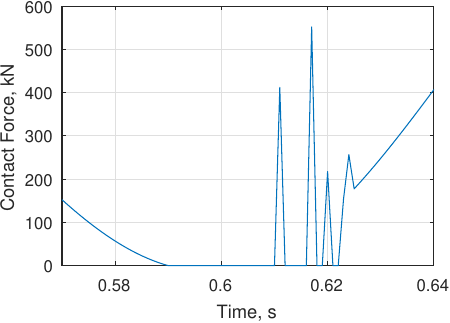}} &
\subfloat[] 
{\label{fig:LCP:detail:displ}\includegraphics{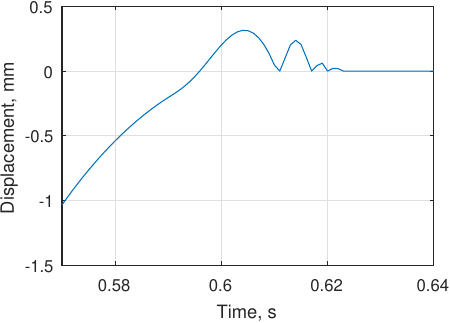}}
\end{tabular}
\caption{Case 6. Contact force $\lambda_1$ (a) and displacements of wheel \#1 (b) during the contact separation}
\label{fig:LCP:detail}
\end{figure}

\noindent ${k^s=50000 \text{kN/m}}$. Mass of the wheel, ${m^w=1000}$kg. Damping of the suspensions is assumed to be equal to $5\%$ of critical damping. Track irregularities are not considered, therefore ${\boldsymbol \rho}(t)={\bf 0}$. Speed of the train ${v_x^{\text t}}$ is constant and equal to ${110 \text{m/s}}$.

Results are presented in Figures \ref{fig:LCP} and \ref{fig:LCP:detail}. It can be seen in Figures \ref{fig:LCP:displ} and \ref{fig:LCP:detail:displ} that wheel \#1 separates from the bridge when leaving the second span of the bridge. At the same time, Figures \ref{fig:LCP:forces} and \ref{fig:LCP:detail:forces} show that contact force ${\lambda_1}$ is equal to zero during the corresponding time interval ${t=0.590...0.622}$s, except the moments when the wheel is bouncing, that is, detaching and reattaching to the track.    

\begin{figure}[H] 
\centering
\captionsetup{justification=centering}
\captionsetup[subfigure]{justification=centering}
\begin{tabular}{l@{\hskip 1.2cm}r}
\subfloat[Contact forces] 
{\label{fig:LCP:Irreg:a}\includegraphics{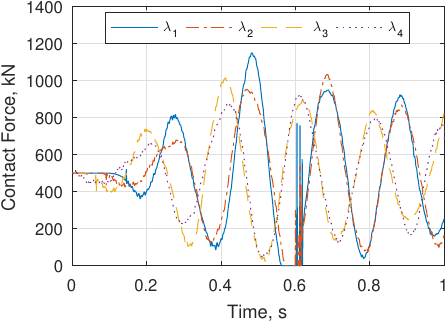}} &
\subfloat[Contact forces $\lambda_1$ and $\lambda_2$ during the contact separation of wheels \#1 and \#2] 
{\label{fig:LCP:Irreg:b}\includegraphics{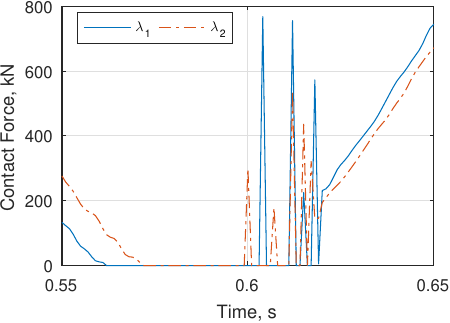}}
\end{tabular}
\caption{Case 6. Amplification of contact forces between the train and the bridge due to the track irregularities}
\label{fig:LCP:Irreg}
\end{figure}

A track irregularity profile (Figure \ref{fig:irreg:b}) is then applied to the bridge. The obtained contact forces are presented in Figure \ref{fig:LCP:Irreg}. As can be seen, not only the last wheel, wheel \#1, separates from the bridge, but wheel \#2 separates as well due to the amplification of contact forces caused by the irregularities of the track.
\FloatBarrier
\section{Preliminary comparison with measured data} \label{sec:compare}

In order to examine the feasibility of the proposed approach, the VTSI algorithm is employed to analyze bridge-train system described by Lee~et~al.~\cite{Lee2010}, for which the measured data is available. The measurements are from field tests conducted by the Korean Institute of Construction Technology \cite{Lee2010}.
In this study, for simplicity, we reduce three-dimensional models to two-dimensional and construct vehicle matrices \Mt, \Ct, \Kt, \Lt\ and \Pt\ accounting for the articulated system of the high-speed train, that is DOF are shared between adjacent passenger cars \cite{Lee2010}. The train is composed of 20 cars, each car has total of 10 DOF. The bridge used in the study is a continuous two-span (40m each) bridge section with a concrete box girder superstructure. The simplified bridge model is assembled with prismatic beam elements as described in Section \ref{sec:2D:bridge:example}. The weight of the bridge is calculated based on the available data. The beam element stiffness is chosen to obtain the first natural frequency of the bridge reported by Lee~et~al.~\cite{Lee2010}.

\begin{figure}[H] 
\centering
\captionsetup{justification=centering}
\captionsetup[subfigure]{justification=centering}
\begin{tabular}{l@{\hskip 1.2cm}r}
\subfloat[Vehicle speed of 300km/h] 
{\label{fig:exp:300}\includegraphics{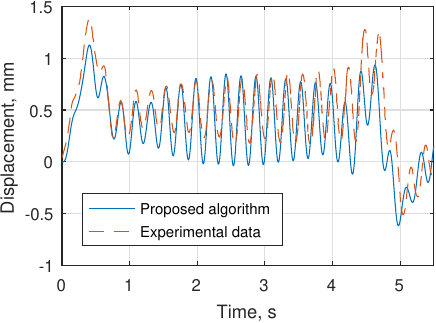}} &
\subfloat[Vehicle speed of 250km/h] 
{\label{fig:exp:250}\includegraphics{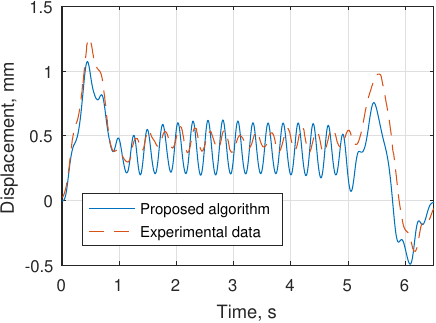}}
\end{tabular}
\caption{Comparison with the experimental measurements (Lee et al. \cite{Lee2010}). Relative vertical displacements of the~middle point of the first span}
\label{fig:exp}
\end{figure}

Figure \ref{fig:exp} illustrates the results of preliminary comparison with experimental data for vehicle speeds of 300km/h and 250km/h. As can be seen, the assembled articulated train model allows to capture the features of bridge displacements, including the peak displacements due to the heaviest cars (the~first and the last two cars) passing the point of measurement. The differences in the obtained results can be attributed to the fact that the bridge model is assembled using two-node prismatic beam elements, the track model is disregarded and other effects like damping in the ballast are ignored. Detailed two- and three-dimensional modeling of the bridge and the track subsystem, as well as detailed experimental validation, are topics of current research. 
\section{Concluding remarks}

An algorithm for two-dimensional vehicle-track-structure interaction (VTSI) analysis has been presented. The equations of motion are derived separately for the train and the bridge subsystems. The coupling of the subsystems is then achieved by enforcing kinematic constraints. The obtained system of equations of motion is discretized in time using self-stabilized discretization scheme proposed by Bauchau \cite{bauchau2003}. Several numerical examples have been analyzed to explore the performance of this algorithm. It is observed that numerical characteristics of response of both, structure and vehicle, subsystems depend on the mass of the wheels. The obtained results also reveal the instability of the numerical scheme in the form of spurious oscillations in the contact forces and accelerations due to curvature discontinuity in the interpolation of kinematic constraints. These oscillations prevent proper evaluation of wheel/track safety and assessment of passenger comfort, as well as validation of the obtain results by experimental data. Therefore, these numerical effects should be eliminated or properly separated from the results of VTSI analysis. Two approaches have been introduced to eliminate these oscillations.

Firstly, a B-spline interpolation of the kinematic constraints has been developed. This interpolation method has been used together with the self-stabilized scheme proposed by Bauchau. Numerical examples demonstrate that this mitigates spurious oscillations, confirming that the reason for their occurrence is curvature discontinuity. However, the B-spline curve does not pass through the nodes of the bridge model, and therefore convergence with increasing number of elements is slow. Moreover, the implementation is complicated and limited by the boundary conditions of the structure. Therefore, the second approach has been developed.

The second approach utilizes a composite time integration scheme proposed by Bathe~\cite{bathe2005}. This two-step scheme remains stable due to the numerical damping introduced at the second sub-step. Numerical examples show that the implemented integration scheme eliminates spurious oscillations in the contact forces and accelerations. When the wheel is entering or leaving the bridge, numerical errors (``undershoots'') in the solutions are observed. These are deemed inconsequential and disregarded. The proposed algorithm has been verified using a generic DAE solver.

Due to high train speeds and possible geometric irregularities in the track geometry, the wheels may ``jump'' off the track. Therefore, a possibility of detachment and subsequent reattachment of the wheels has been incorporated into the algorithm. A Linear Complementary Problem (LCP) solver has been employed to modify the procedure.

Modeling the bridge and the train as independent objects makes the proposed algorithm highly modular and allows its incorporation into existing bridge analysis software without significant modification. As a result, a variety of tools in the existing software can be used to create a bridge model, e.g.  truss, arch, beam or cable stayed bridge, single- or multi-span, fixed or simply supported, path definitions, influence surface generation. The required train model can be created in any multibody dynamics software  or  specialized  railway  vehicle  design  software  and  imported  into  bridge  analysis software in form of structural matrices. Although only simple bridge models have been considered in the examples presented in this paper, when coupled with fully featured bridge analysis software, complex and diverse bridge models can be used without any change to the algorithm.

The main contribution of this paper is recognizing computational issues associated with time-varying kinematic constraints with discontinuous derivatives, clearly identifying their cause and developing remedies. While the focus of this paper is on the computational procedure, detailed modeling of the bridge and track, modeling curved paths and rolling conditions, and practical applications of VTSI analysis are topics of ongoing research.

\appendix





\bibliographystyle{elsarticle-num}

\bibliography{references}

\begin{thebibliography}{10}
\expandafter\ifx\csname url\endcsname\relax
  \def\url#1{\texttt{#1}}\fi
\expandafter\ifx\csname urlprefix\endcsname\relax\def\urlprefix{URL }\fi
\expandafter\ifx\csname href\endcsname\relax
  \def\href#1#2{#2} \def\path#1{#1}\fi

\bibitem{TM}
\textit{Parsons Brinckerhoff}, Track-structure interaction, Technical Memorandum TM 2.10.10, California High-Speed Train Project (2013).

\bibitem{neves2012direct}
S.~Neves, A.~Azevedo, R.~Calcada, A direct method for analyzing the vertical vehicle-structure interaction, Engineering Structures 34 (2012) 414--420.

\bibitem{Montenegro2016}
P.~A. Montenegro, R.~Cal{\c{c}}ada, N.~{Vila Pouca}, M.~Tanabe, {Running safety assessment of trains moving over bridges subjected to moderate earthquakes}, Earthquake Engineering and Structural Dynamics 45 (2016) 483--504.
\newblock \href {https://doi.org/10.1002/eqe.2673} {\path{doi:10.1002/eqe.2673}}.

\bibitem{Xia2014}
C.~Y. Xia, H.~Xia, G.~{De Roeck}, {Dynamic response of a train-bridge system under collision loads and running safety evaluation of high-speed trains}, Computers and Structures 140 (2014) 23--38.

\bibitem{Sun}
Z.~Sun, Z.~Zou, {Towards an efficient method of predicting vehicle-induced response of bridge}, Engineering Computations 33~(7) (2016) 2067--2089.

\bibitem{Xia2005}
H.~Xia, N.~Zhang, {Dynamic analysis of railway bridge under high-speed trains}, Computers {\&} Structures 83 (2005) 1891--1901.

\bibitem{Dimitrakopoulos2015}
E.~G. Dimitrakopoulos, Q.~Zeng, {A three-dimensional dynamic analysis scheme for the interaction between trains and curved railway bridges}, Computers {\&} Structures 149 (2015) 43--60.

\bibitem{Salcher2015}
P.~Salcher, C.~Adam, {Modeling of dynamic train-bridge interaction in high-speed railways}, Acta Mechanica 226~(8) (2015) 2473--2495.

\bibitem{Yang1999}
Y.-B. Yang, C.-H. Chang, J.-D. Yau, {An element for analysing vehicle-bridge systems considering vehicle's pitching effect}, International Journal for Numerical Methods in Engineering 46 (1999) 1031--1047.

\bibitem{Diana1989}
G.~Diana, F.~Cheli, {Dynamic Interaction of Railway Systems with Large Bridges}, Vehicle System Dynamics 18 (1989) 71--106.

\bibitem{Kwark2004}
J.~W. Kwark, E.~S. Choi, Y.~J. Kim, B.~S. Kim, S.~I. Kim, {Dynamic behavior of two-span continuous concrete bridges under moving high-speed train}, Computers {\&} Structures 82~(4-5) (2004) 463--474.

\bibitem{Zhang2013}
N.~Zhang, H.~Xia, {Dynamic analysis of coupled vehicle-bridge system based on inter-system iteration method}, Computers {\&} Structures 114-115 (2013) 26--34.

\bibitem{sivaselvan2014}
M.~Sivaselvan, J.~Tauberer, A.~Karakaplan, Dynamic vehicle-track-structure interaction analysis using lagrange multipliers, in: Proceedings of the Istanbul Bridge Conference, 2014.

\bibitem{Ferris1997}
M.~C. Ferris, J.~S. Pang, {Engineering and Economic Applications of Complementarity Problems}, SIAM Review 39~(4) (1997) 669--713.

\bibitem{Kwak1991}
B.~M. Kwak, {Complementarity Problem Formulation of Three-Dimensional Frictional Contact}, Journal of Applied Mechanics 58 (1991) 134--140.

\bibitem{Xi2016_2D}
Y.~Xi, A.~Almqvist, Y.~Shi, J.~Mao, R.~Larsson, {A Complementarity Problem-Based Solution Procedure for 2D Steady-State Rolling Contacts with Dry Friction}, Tribology Transactions 59~(6) (2016) 1031--1038.

\bibitem{Xi2016_3D}
Y.~Xi, A.~Almqvist, Y.~Shi, J.~Mao, R.~Larsson, {Linear Complementarity Framework for 3D Steady-State Rolling Contact Problems Including Creepages with Isotropic and Anisotropic Friction for Circular Hertzian Contact}, Tribology Transactions (2016) 1--13\href {https://doi.org/10.1080/10402004.2016.1217111} {\path{doi:10.1080/10402004.2016.1217111}}.

\bibitem{acary2008numerical}
V.~Acary, B.~Brogliato, Numerical Methods for Nonsmooth Dynamical Systems: Applications in Mechanics and Electronics, Lecture Notes in Applied and Computational Mechanics, Springer Berlin Heidelberg, 2008.

\bibitem{Zhu2015}
D.~Y. Zhu, Y.~H. Zhang, H.~Ouyang, {A linear complementarity method for dynamic analysis of bridges under moving vehicles considering separation and surface roughness}, Computers {\&} Structures 154 (2015) 135--144.

\bibitem{brenan1996numerical}
K.~E. Brenan, S.~L.~V. Campbell, S.~L. Campbell, L.~R. Petzold, {Numerical Solution of Initial-Value Problems in Differential-Algebraic Equations}, Classics in Applied Mathematics, Society for Industrial and Applied Mathematics, 1996.

\bibitem{bauchau2003}
O.~Bauchau, A self-stabilized algorithm for enforcing constraints in multibody systems, International journal of solids and structures 40~(13) (2003) 3253--3271.

\bibitem{yang2004vehicle}
Y.-B. Yang, J.~Yau, Y.~Wu, Vehicle-bridge interaction dynamics with applications to high-speed raileawy, World Scientific, 2004.

\bibitem{spiryagin2014design}
M.~Spiryagin, C.~Cole, Y.~Sun, M.~McClanachan, V.~Spiryagin, T.~McSweeney, Design and simulation of rail vehicles, CRC Press, 2014.

\bibitem{rogers2000introduction}
D.~Rogers, An introduction to NURBS: with historical perspective, Elsevier, 2000.

\bibitem{bathe2005}
K.-J. Bathe, M.~Baig, On a composite implicit time integration procedure for nonlinear dynamics, Computers \& Structures 83~(31) (2005) 2513--2524.

\bibitem{bathe2012}
K.-J. Bathe, G.~Noh, Insight into an implicit time integration scheme for structural dynamics, Computers \& Structures 98 (2012) 1--6.

\bibitem{serban2015user}
R.~Serban, C.~Petra, A.~Hindmarsh, User Documentation for idas v1. 2.2 (sundials v2. 6.2), Center for Applied Scientific Computing. Lawrence Livermore National Laboratory, 2015.

\bibitem{SUNDIALS}
R.~Serban, \href{http://computation.llnl.gov/projects/sundials-suite-nonlinear-differential-algebraic-equation-solvers/sundials-software}{A matlab interface to sundials} (2015).
\newline\urlprefix\url{http://computation.llnl.gov/projects/sundials-suite-nonlinear-differential-algebraic-equation-solvers/sundials-software}

\bibitem{Claus1998}
H.~Claus, W.~Schiehlen, {Modeling and simulation of railway bogie structural vibrations}, Vehicle System Dynamics 29 (1998) 538--552.

\bibitem{cottle1992linear}
R.~W. Cottle, J.~S. Pang, R.~E. Stone, {The Linear Complementarity Problem}, Classics in Applied Mathematics, Society for Industrial and Applied Mathematics, 1992.

\bibitem{LCP}
A.~Almqvist, A.~Spencer, P.~Wall, \href{http://www.mathworks.com/matlabcentral/fileexchange/41485-a-pivoting-algorithm-solving-linear-complementarity-problems/content/LCPSolve.m}{A pivoting algorithm solving linear complementarity problems} (2013).
\newline\urlprefix\url{http://www.mathworks.com/matlabcentral/fileexchange/41485-a-pivoting-algorithm-solving-linear-complementarity-problems/content/LCPSolve.m}

\bibitem{Lee2010}
Y.-S. Lee, S.-H. Kim, {Structural analysis of 3D high-speed train-bridge interactions for simple train load models}, Vehicle System Dynamics 48~(2) (2010) 263--281.

\end{thebibliography}

\end{document}